\documentclass[11pt,leqno]{amsart}
\usepackage{graphicx}
\usepackage{indentfirst,csquotes}

\usepackage{amssymb,amsthm,amsmath}
\usepackage{xcolor,hyperref,fancyhdr,etoolbox}
\usepackage{booktabs,enumitem,algorithm,algpseudocode}

\newtheorem{theorem}{Theorem}
\newtheorem{fact}{Fact}
\newtheorem{lemma}{Lemma}
\theoremstyle{definition}
\newtheorem{definition}{Definition}
\theoremstyle{remark}
\newtheorem{remark}{Remark}

\hypersetup{ colorlinks=true, linkcolor=black, filecolor=black, urlcolor=black, citecolor=black }

\newcommand{\citep}[1]{\cite{#1}}

\begin{document}
\title{Gradient Discrete Schr\"odinger Bridge Matching for Peptide Generation}
\author[A. BenAbdallah, M. Dlala]{Abdallah BenAbdallah}
\address{Department of Mathematics, Higher Institute of Computer Science and Multimedia of Sfax, University of Sfax, Tunisia}
\email{abdellah.benabdallah@isims.usf.tn}
\author{Mohsen Dlala}
\address{Department of Mathematics, Qassim University, Buraydah, Saudi Arabia}
\email{3862@qu.edu.sa}
\date{\today}
\maketitle

{\let\thefootnote\relax}

\begin{abstract}
The optimal control of a discrete Schr\"odinger bridge is a discrete gradient, a
difference of vertex potentials, but recent models such as MadSBM instead learn a
free control on the edges by cross-entropy. We prove that this leaves the action
uncontrolled: cross-entropy identifies the control only up to a state-dependent
shift (a gauge), and the action is $\operatorname{KL}(r\|p)$ plus a gauge term
the loss never constrains. We propose GDSBM, which writes the control as a
difference of vertex potentials and adds the action to the loss; pinning the
gauge costs no expressiveness and makes the action a property of the model, so
the penalty acts on the model and not on the gauge. On peptide generation our
MadSBM reimplementation has a free-gauge action of order $10^{8}$ for a
gauge-fixed value of $6.5$, against $1.5$ for GDSBM, with no quality loss that
three seeds can detect. The gradient form also fits guidance, which becomes the
addition of a potential: guided toward a net charge of $+3$ at pH~7, GDSBM is
closer than MadSBM, over three seeds, to its exact tempered target in total
variation at every guidance strength (at $\beta=1$, $81\%$ of samples in the
target window against $74\%$, for a target of $97\%$); we claim no statistical
significance.
\end{abstract}

\bigskip

\section{Introduction}\label{sec:intro}

Peptides occupy the therapeutic space between small molecules and antibodies:
specific enough to hit flat protein surfaces, small enough to manufacture
chemically \citep{muttenthaler2021}. Their discovery still relies on screening
large libraries and keeping the rare hits, whereas generative models promise to
invert the pipeline and sample candidates \emph{conditioned on the desired
function} \citep{alamdari2023evodiff,watson2023rfdiffusion,tang2024peptune}. The
promise is already concrete for proteins---structure prediction
\citep{jumper2021highly,abramson2024accurate}, language models
\citep{rives2021biological,lin2023esm,madani2023large,hayes2025simulating},
structure- and sequence-level design
\citep{dauparas2022robust,ingraham2023illuminating}---and for antimicrobial
peptides, where generative models have produced laboratory-tested candidates
\citep{das2020accelerating,szymczak2023discovering} and biologically informed
generators are emerging \citep{soares2025omegamp}. Making it quantitative
requires two things current discrete generators do not simultaneously offer: a
principled training objective, and conditioning that is exact rather than
heuristic.

Continuous-state models dominate protein \emph{structure} generation
\citep{ho2020denoising,song2020score,liu2022flow,albergo2022building,yim2023se,lin2023generating,bose2024se,watson2023rfdiffusion};
on \emph{sequences} the generators of choice are discrete diffusions and flows
\citep{austin2021d3pm,hoogeboom2021argmax,campbell2022,lou2024sedd,sahoo2024mdlm,shi2024simplified,campbell2024multiflow,stark2024dirichlet},
now at language-model scale \citep{nie2026large} and instantiated for proteins
\citep{alamdari2023evodiff,wang2024diffusion,lisanza2025multistate}. Closest to
this work is Schr\"odinger-bridge matching, first continuous
\citep{debortoli2021,shi2023,tong2024} and now discrete
\citep{kim2025discrete}---but validated only on image tokens and graphs
\citep{kim2025discrete,ksenofontov2025categorical}, with closed-form benchmarks
only now appearing \citep{aramayo2026entering}.

\begin{figure}[t]
\centering
\includegraphics[width=0.8\linewidth]{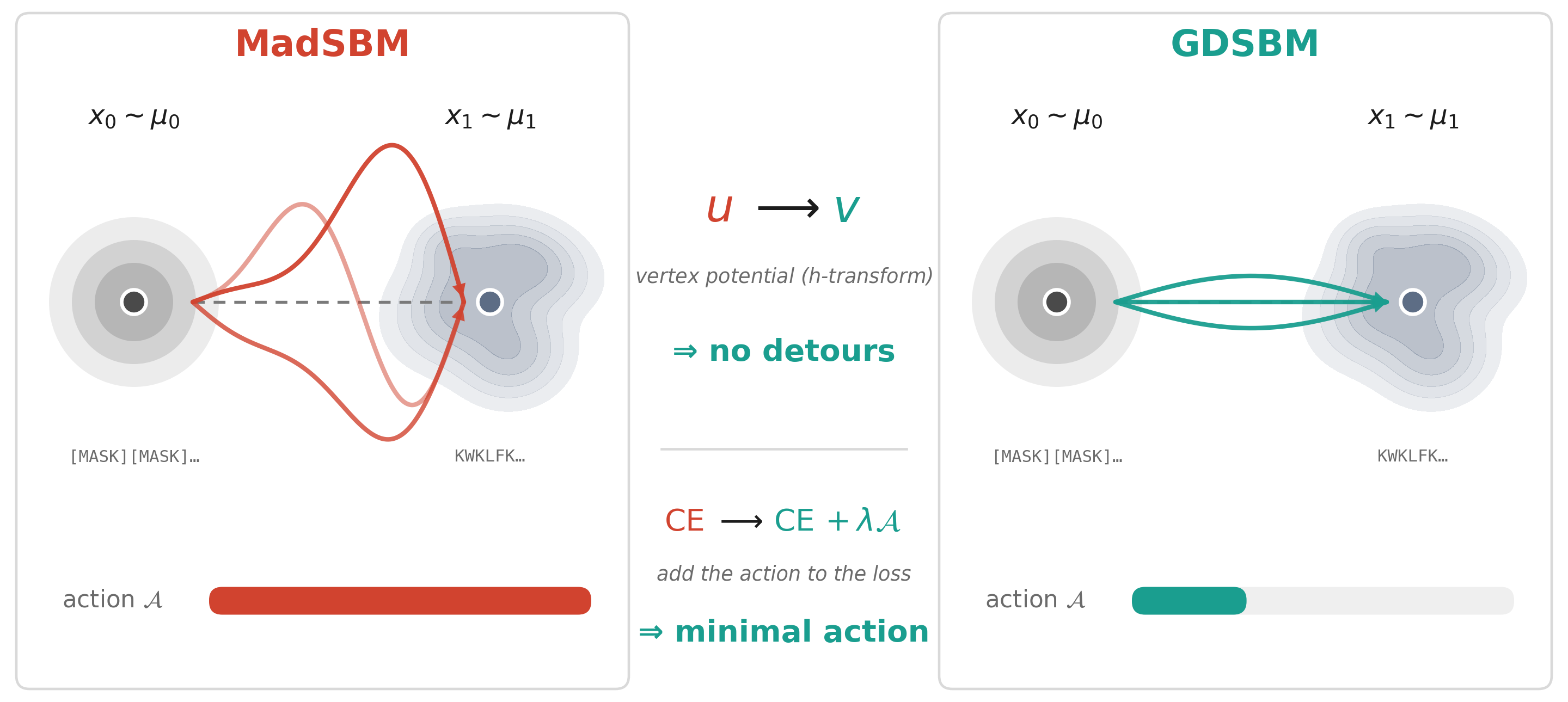}
\caption{\textbf{Two changes turn MadSBM into a bridge.} Both models transport a fully
masked $x_0$ to a peptide $x_1$; bar lengths are illustrative (measured values in
Section~\ref{sec:experiments}). \emph{Left, MadSBM:} $u$ is a free field on the edges,
trained by cross-entropy, which fixes it only up to a state-dependent shift and leaves the
action free. \emph{Right, GDSBM:} (i)~$u=\bar\nabla v$, the form of the Doob optimum, which
has zero circulation along cycles; (ii)~$\mathrm{CE}\to\mathrm{CE}+\lambda\mathcal A$, which
lowers the gauge-fixed action ($6.35$ against $1.51$ at $N=32$).}
\label{fig:teaser}
\end{figure}

MadSBM \citep{goel2026madsbm} is, to our knowledge, the first discrete bridge
model brought to the scale of biological sequences: a controlled jump process on
the amino-acid edit graph, anchored to a frozen protein language model
\citep{lin2023esm}, with the forward--backward iterations
\citep{theodoropoulos2026generalized} that are hard to scale to
$|\mathcal V|^L$ states replaced by a token-wise cross-entropy. This is a
genuine advance and our starting point. Its price, however, is structural, and
it is precisely what makes a bridge a bridge. \emph{The control is not an
$h$-transform of the ESM-2 reference}: it is a free field $u_\theta(x,y,t)$ on
the \emph{edges} which, by MadSBM's own consistency theorem, cross-entropy
recovers only up to an additive $c(x,t)$, whereas the optimum is a difference of
\emph{vertex} values, $u^*(x,y,t)=\log h(y,t)-\log h(x,t)$. \emph{The loss is not
the action}: cross-entropy is minimized by every representative of the class
$u+c$, whose actions differ without bound (Section~\ref{sec:gauge}), so the
reported action is a statement about the representative the optimizer happens to
produce rather than about the model. Conditioning, finally, is a logit tilt known
to be biased for tempered targets \citep{lee2025debiasing,wan2026dgm}, a bias
usually asserted and left unmeasured. MadSBM keeps the construction and the name
of the Schr\"odinger bridge but, relative to the ESM-2 reference, neither its
control class nor its objective. Figure~\ref{fig:teaser} contrasts the three
objects.

\paragraph{Observation, and a gauge.} Doob's $h$-transform states that the
optimal bridge control is a difference of \emph{vertex} functions
\citep{follmer2006random,leonard2014}. We take this seriously three times: for
\emph{parameterization} (learn a vertex potential, not an edge field), for
\emph{conditioning} (design objectives are vertex functions, and composing
$h$-transforms adds potentials), and for \emph{inference-time correction} (the
twisting functions of sequential Monte Carlo are vertex functions). A fourth
observation concerns not what the control should be but what training
\emph{identifies}. At each state the conditional law of the next token is a
softmax of reference logits plus control, and a softmax is invariant to adding a
function of the current state; a control learned by cross-entropy is therefore
identified only up to this \emph{gauge}, whereas the action is not invariant to
it. It splits exactly into $\operatorname{KL}(r\|p_\theta)$, a function of the
learned conditional alone, and an unbounded gauge excess
(Theorem~\ref{thm:gauge}). The vertex form of the optimum is the natural gauge
fixing: it pins the term the softmax cannot see, so that an action penalty acts
on the conditional and not on the gauge. This changes how the action of a bridge
model must be reported, and what it means to regularize it.

\paragraph{Contributions.} Three main results, one per section, with all proofs in the appendix.
\begin{itemize}[leftmargin=1.4em,itemsep=0.1em,topsep=0.2em]
\item \textbf{The action of a free control measures its gauge (C1,
  Theorem~\ref{thm:gauge}).} The per-state action splits exactly into
  $\operatorname{KL}(r\|p_\theta)$, a function of the learned conditional alone,
  and a gauge excess that cross-entropy never constrains and that is unbounded
  along the orbit of a model; at the forward Girsanov kernel the gauge-fixed
  value is instead bounded and invertible in $\operatorname{KL}(p_\theta\|r)$.
  Measured, a trained free field reports $10^{8}$ for a gauge-fixed $6.5$.
\item \textbf{Vertex-potential control and the exact effect of the penalty (C2,
  Theorem~\ref{thm:vertex}).} We parameterize
  $u_\theta(x,y,t)=v_\theta(y,t)-v_\theta(x,t)$, the form of the optimum, from a
  zero-initialized energy head. This pins the gauge without restricting the
  representable conditionals, so at $\lambda=0$ the model is functionally
  identical to MadSBM \citep{goel2026madsbm}; the gauge term becomes an anchor on
  the current token; and the penalty has an exact effect, a convex shrinkage of
  the data posterior towards the reference with weight
  $\lambda\nu_V/(1+\lambda\nu_V)$.
\item \textbf{Conditioning is algebra on potentials, and its bias is measurable
  (C3, Theorem~\ref{thm:cond}).} A tilt is an $h$-transform whose potential adds
  to the trained one, so a multi-objective sweep needs no retraining; for any
  scalar statistic the tempered target's marginal is closed-form, with no product
  structure assumed; and a particle system with telescoping weights reaches it.
  On peptides the net charge is exactly computable, turning the bias of the
  myopic guided sampler into a measurement against an explicit sampling floor.
\end{itemize}

\paragraph{What we do not claim.} No wet-lab validation and no state-of-the-art
claim: we compare to our own MadSBM reimplementation at matched compute and
seeds, reporting effect sizes and sign consistency, not significance. We do not
claim that $v_\theta$ is a globally consistent potential (consistency across
states is learned, not enforced), that the vertex form helps by itself (at
$\lambda=0$ the two models coincide) or beyond the penalty (the $2\times2$
ablation has one seed, $300$ steps), or that we minimize the Schr\"odinger
action. Multi-objective guidance and the particle system are stated, not evaluated.

\paragraph{Outline.} Section~\ref{sec:related} situates related work;
Section~\ref{sec:background} recalls bridges on graphs;
Sections~\ref{sec:gauge}--\ref{sec:conditioning} state the three results;
Section~\ref{sec:experiments} reports the action audit, the form/penalty
ablation, unconditional generation and exactly-conditioned peptide generation.
Deferred remarks are collected in Appendix~\ref{app:remarks}.

\section{Related Work}\label{sec:related}
Our question is narrow: for sequence generation on a discrete state space, should
the control of a Schr\"odinger bridge be a free field on the edges of the sequence
graph, or an exact discrete gradient of a vertex potential?

\paragraph{Bridges, and potentials versus fields.} The Schr\"odinger problem
\citep{schrodinger1932,leonard2014} is the entropic regularisation of dynamical
optimal transport \citep{benamou2000computational,chen2021stochastic,cuturi2013};
continuous-state neural solvers differ in which fixed-point iteration they
discretise \citep{debortoli2021,shi2023,tong2024,liu2024generalized}. In
$\mathbb R^d$ the optimal control is a gradient, $u^\ast=\nabla\log h$, and
methods that exploit this---scalar-potential transport
\citep{neklyudov2023action}, adjoint and stochastic-optimal-control matching
\citep{domingo2025adjoint,domingo2026adjoint}, energy-based parametrisations
\citep{balcerak2026energy}---gain by removing a rotational component that carries
no transport. On graphs the geometry of entropy-driven dynamics is well
understood \citep{maas2011gradient,erbar2012ricci,chow2012fokker} and is only now
becoming a learning algorithm \citep{rancati2026learning}, but on a simplex with
$(|\mathcal V|+1)^L$ vertices neither the metric nor the free energy is
computable, so we work with the control of a jump process and the discrete
gradient $v(y)-v(x)$ along an edge.

\paragraph{Discrete generation and discrete bridges.} Discrete diffusion
\citep{austin2021d3pm,campbell2022}, score-entropy models \citep{lou2024sedd},
masked diffusion \citep{sahoo2024mdlm,shi2024simplified} and discrete flow
matching \citep{campbell2024multiflow} learn the reversal of a reference jump
process; discrete bridge formulations
\citep{kim2025discrete,ksenofontov2025categorical,aramayo2026entering,theodoropoulos2026generalized}
place a \emph{control} on its rates. MadSBM \citep{goel2026madsbm} does so for
masked sequences with an ESM-2 reference \citep{lin2023esm}, trained by
cross-entropy; it is our direct baseline. All learn general rate structures,
whereas we impose $u_\theta(x,y,t)=v_\theta(y,t)-v_\theta(x,t)$, which removes
cycle inconsistencies by construction, and we make the associated action explicit
and use it to audit the learned bridge.

\paragraph{Guidance, reward optimisation, and the test domain.}
Function-conditioned design has been attacked as model-based optimisation and
reinforcement learning
\citep{brookes2019cbas,angermueller2020dynappo,stanton2022accelerating,sinai2020adalead,kirjner2024improving,bengio2021flow,jain2022biological,notin2023proteingym},
as inference-time guidance
\citep{dhariwal2021diffusion,ho2022classifier,chung2022diffusion,gruver2023protein,nisonoff2024,schiff2024,rectorbrooks2024,li2024derivative,wu2023tds,trippe2023smcdiff,skreta2025feynman},
and as reward fine-tuning
\citep{black2024training,uehara2024understanding,domingo2025adjoint,wang2025fine},
all presuming the generator is a well-defined object whose conditional laws can be
tilted; whether the \emph{control} of a bridge is such an object is our question.
The exact twist, $h(x,t)=\mathbb E[e^{\varphi(X_T)}\mid X_t=x]$, is intractable
here and practical surrogates are biased \citep{lee2025debiasing}; we apply the
same $h$-transform viewpoint to a control that is already a vertex potential, so
the twist adds a second potential to the first. We evaluate on peptides, as MadSBM does, for their natural reference process and
exactly computable objectives; our comparison is to MadSBM alone, at matched
corpus, compute and seeds.

\section{Background: Bridges on the Sequence Graph}\label{sec:background}

This section fixes the objects used throughout and recalls two classical facts:
the Girsanov formula for jump processes (Fact~\ref{fact:girsanov}) and the
Gibbs--Doob variational principle, which makes optimal controls \emph{gradients}
(Fact~\ref{fact:gibbs}). Neither is a contribution; both are proved in
Appendices~\ref{app:girsanov} and~\ref{app:gibbs} because the paper reads them
closely.

\paragraph{Sequences as a jump process.} Let $V$ be the set of $20$ amino acids,
$\mathsf m\notin V$ a mask token, and $\mathcal X=(V\cup\{\mathsf m\})^L$. The
\emph{edit graph} $\mathcal E\subset\mathcal X\times\mathcal X$ contains the pairs
$(x,y)$ such that $y$ is obtained from $x$ by replacing one masked position $i$
by one amino acid $a$; we write $y=x^{i\leftarrow a}$. The reference process is
the continuous-time Markov chain on $\mathcal X$ with rates
\begin{equation}\label{eq:ref-rates}
R^0_t(x,y)=\lambda_t\,r_t(a\mid x,i)\,\mathbf 1_{\{x_i=\mathsf m\}},\qquad y=x^{i\leftarrow a},
\end{equation}
and $R^0_t(x,y)=0$ off $\mathcal E$, where $\lambda_t\ge0$ is the masking
intensity and $r_t(\cdot\mid x,i)$ a probability vector on $V$ (built from the
protein language model in Section~\ref{sec:gdsbm}). We fix a horizon $T$ on which
the rates are bounded, and write $\mathbb P^0$ for the law on
$D([0,T],\mathcal X)$ of the chain started from the Dirac mass at $\mathsf m^L$.
Since $\mathcal X$ is finite and the rates are bounded, all objects below are
well defined.

\paragraph{Controls and action.} A \emph{control} is a bounded measurable
$u:[0,T]\times\mathcal E\to\mathbb R$; the controlled chain has rates
$R^u_t(x,y)=R^0_t(x,y)e^{u_t(x,y)}$ and the same initial law, with path law
$\mathbb P^u$. The cost of a control is the relative entropy of its path law.

\begin{fact}[Girsanov for jump processes]\label{fact:girsanov}
Let $(t_k,X_{t_k})_{k\le n}$ be the jump times and states of a path. Then $\mathbb P^u\ll\mathbb P^0$ with
\begin{equation}\label{eq:girsanov-density}
\begin{split}
\log\frac{d\mathbb P^u}{d\mathbb P^0}={}&\sum_{k=1}^{n}u_{t_k}(X_{t_k-},X_{t_k})\\
&-\int_0^T\sum_y R^0_t(X_t,y)\big(e^{u_t(X_t,y)}-1\big)\,dt,
\end{split}
\end{equation}
and, with $a(z)=ze^z-e^z+1$ and $\Psi(z)=e^z-z-1$,
\begin{align}
\operatorname{KL}(\mathbb P^u\,\|\,\mathbb P^0)&=\mathbb E_{\mathbb P^u}\!\int_0^T\sum_y R^0_t(X_t,y)\,a\big(u_t(X_t,y)\big)\,dt,\label{eq:kl-forward}\\
\operatorname{KL}(\mathbb P^0\,\|\,\mathbb P^u)&=\mathbb E_{\mathbb P^0}\!\int_0^T\sum_y R^0_t(X_t,y)\,\Psi\big(u_t(X_t,y)\big)\,dt.\label{eq:kl-reverse}
\end{align}
\end{fact}

Both identities are elementary once the path likelihood of a finite-state jump
process is written down: \eqref{eq:girsanov-density} is a ratio of two such
likelihoods, and \eqref{eq:kl-forward}--\eqref{eq:kl-reverse} follow by taking
expectations and compensating the jump term, once under $\mathbb P^u$ and once
under $\mathbb P^0$. It is this change of compensating measure, and nothing else,
that produces the two kernels.

\begin{remark}[Two kernels]\label{rem:two-kernels}
The kernels are related by $\Psi(z)=e^{z}a(-z)$. Both are non-negative and vanish
only at $z=0$, but they are \emph{not} interchangeable: $\Psi$ is strictly
convex, while $a''(z)=(z+1)e^z$ is negative for $z<-1$, so $a$ is not convex.
Equation~\eqref{eq:kl-forward}, with $a$, is the path-space action
$\operatorname{KL}(\mathbb P^u\|\mathbb P^0)$ that defines the Schr\"odinger
bridge; \eqref{eq:kl-reverse}, with $\Psi$, is its reverse
counterpart.\footnote{The published version of MadSBM writes $ze^z-e^z-1$
(Corollary~4.2); the identity of their Theorem~4.1 gives $+1$, which is required
for the kernel to vanish at $z=0$. We use $+1$.} Section~\ref{sec:gdsbm} states
which kernel the training penalty uses and in which expectation; we do not
identify a per-state penalty with a path-space entropy.
\end{remark}

\paragraph{Doob's transform.} Let $\varphi:\mathcal X\to\mathbb R$ be a terminal
reward and set $h(x,t)=\mathbb E_{\mathbb P^0}[e^{\varphi(X_T)}\mid X_t=x]$,
$v_t(x)=\log h(x,t)$. Let $\mathbb P^h$ be the chain with rates
$R^0_t(x,y)h(y,t)/h(x,t)$, i.e.\ the controlled chain with
\begin{equation}\label{eq:doob-control}
u^\ast_t(x,y)=v_t(y)-v_t(x), \qquad (x,y)\in\mathcal E.
\end{equation}

\begin{fact}[Gibbs--Doob variational principle]\label{fact:gibbs}
$\dfrac{d\mathbb P^h}{d\mathbb P^0}=\dfrac{e^{\varphi(X_T)}}{h(x_0,0)}$, and for every path law $\mathbb P\ll\mathbb P^0$ with the same initial state $x_0$,
\begin{equation}\label{eq:gibbs}
\operatorname{KL}(\mathbb P\,\|\,\mathbb P^h)=\operatorname{KL}(\mathbb P\,\|\,\mathbb P^0)-\mathbb E_{\mathbb P}\,\varphi(X_T)+\log h(x_0,0)\;\ge 0 .
\end{equation}
Hence $\min_{\mathbb P}\{\operatorname{KL}(\mathbb P\|\mathbb P^0)-\mathbb E_{\mathbb P}\varphi(X_T)\}=-\log h(x_0,0)$, attained uniquely at $\mathbb P^h$.
\end{fact}

The proof is a single cancellation: writing \eqref{eq:girsanov-density} for
$u^\ast$, the jump term telescopes against the increments of $\log h(X_t,t)$, and
the backward equation for $h$ makes what remains exactly the compensator with the
opposite sign (Appendix~\ref{app:gibbs}). For a prescribed terminal law
$\mu_1\ll\mathbb P^0_T$ and a Dirac source the Schr\"odinger system reduces to a
single potential, $\mathbb P^\ast=\frac{d\mu_1}{d\mathbb P^0_T}(X_T)\,\mathbb P^0$
\citep{csiszar1975divergence,leonard2014}. What matters here is the structure of
\eqref{eq:doob-control}: the optimal control is the \emph{discrete gradient} of a
scalar function on the vertices, not an arbitrary function on the
$\mathcal O(L|V|\,|V|^{L})$ edges. It is intractable to compute, since $h$
requires an expectation over $|V|^L$ continuations; GDSBM learns $v$ instead.

\paragraph{Learning the rates: MadSBM and the token kernel.} MadSBM
\citep{goel2026madsbm} instantiates $R^0$ with ESM-2 \citep{lin2023esm} as
reference token distribution $r$ and learns an edge-wise control
$u_\theta(x,y,t)$, trained by cross-entropy on the clean token under the model
posterior
\begin{equation}\label{eq:softmax-post}
p_\theta(y\mid x,t)=\frac{R^0_t(x,y)e^{u_\theta(x,y,t)}}{Z_\theta(x,t)},\qquad Z_\theta(x,t)=\sum_{y'}R^0_t(x,y')e^{u_\theta(x,y',t)}.
\end{equation}
Since \eqref{eq:softmax-post} is invariant under
$u_\theta\mapsto u_\theta+c(x,t)$ for any $c$ independent of $y$, cross-entropy
determines $u$ only modulo such a gauge; their consistency theorem states exactly
this, $u_{\theta^\ast}=u^\ast+c(x,t)$ almost everywhere on $\rho^\ast_t$, under an
expressivity assumption. The \emph{action} \eqref{eq:kl-reverse}, however, is
\emph{not} gauge-invariant.

\paragraph{Per-state notation.} Fix a state $x$ and a time $t$ and set
$\nu=\sum_{y}R^0_t(x,y)$, $r(y)=R^0_t(x,y)/\nu$,
$Z=\sum_{y}R^0_t(x,y)e^{u_t(x,y)}$, $p(y)=R^0_t(x,y)e^{u_t(x,y)}/Z$ and
$\zeta=Z/\nu=\sum_{y}r(y)e^{u_t(x,y)}$. Thus $\nu$ is the total reference rate
out of $x$ and $r$ the reference law of the destination; $\zeta$ is the factor by
which the control inflates the total rate, and $p$ is the law of the destination
of the next jump out of $x$ under the controlled chain, i.e.\ the token kernel
the sampler draws from. Finally $h(\zeta)=\zeta-1-\log\zeta$ for $\zeta>0$ is
non-negative and vanishes only at $\zeta=1$.

\section{The Action of a Cross-Entropy-Trained Control}\label{sec:gauge}

The action \eqref{eq:kl-reverse} is a sum over states of a per-state functional.
The following result computes it in closed form and separates the part the model
determines from the part it does not.

\begin{theorem}[Gauge decomposition of the action]\label{thm:gauge}
With the notation above:
\begin{enumerate}[label=(\alph*),leftmargin=1.8em,itemsep=0.1em,topsep=0.2em]
\item \emph{(splitting)}
\begin{equation}\label{eq:gauge-decomp}
\sum_{y}R^0(x,y)\,\Psi\big(u(y)\big)=\nu\operatorname{KL}(r\,\|\,p)+\nu\,h(\zeta),
\end{equation}
both terms being non-negative; the first depends on $u$ only through $p$, the second only through $\zeta$.
\item \emph{(non-identifiability)} $u\mapsto u+c\mathbf 1$ leaves $p$ (hence the cross-entropy) unchanged and sends $\zeta\mapsto\zeta e^{c}$. Along this orbit the action equals $\nu\operatorname{KL}(r\|p)+\nu\,h(\zeta e^{c})$: it takes every value in $[\nu\operatorname{KL}(r\|p),\infty)$, its minimum being attained at $c^\ast=-\log\zeta$.
\item \emph{(forward kernel)} for $u=\log(p/r)+c\mathbf 1$,
\[
\sum_yR^0(x,y)\,a\big(u(y)\big)=\nu\Big[1+e^{c}\big(\operatorname{KL}(p\|r)+c-1\big)\Big],
\]
minimised at $c^\ast=-\operatorname{KL}(p\|r)$ with value $\nu\big(1-e^{-\operatorname{KL}(p\|r)}\big)\in[0,\nu)$; the map $\operatorname{KL}(p\|r)\mapsto$ this value is a bijection from $[0,\infty)$ onto $[0,\nu)$.
\end{enumerate}
\end{theorem}

The proof is in Appendix~\ref{app:gauge}. Writing $u=\log(p/R^0)+\log Z$ splits
$\Psi(u)=e^{u}-u-1$ into a relative entropy, blind to a constant shift of $u$,
and a term $h(\zeta)$ that sees nothing else: the separation is the softmax
invariance of $p$, read on the action. The number a model reports as its action
is therefore a genuine property, $\nu\operatorname{KL}(r\|p)$, plus a term that
by~(b) can be made arbitrarily large without altering the token kernel at any
state, so any comparison of raw actions is void, including one in our own favour
(Remark~\ref{rem:gauge}).

\section{GDSBM: Vertex-Potential Control with Action Regularization}\label{sec:gdsbm}

Theorem~\ref{thm:gauge} says that the action of a free control is not a model
quantity. This section fixes that by reading the control as a difference of a
vertex potential, the form of the Doob optimum \eqref{eq:doob-control}.

\paragraph{The trained object.} Let $\mathcal T$ be the token vocabulary of the
protein language model (the $20$ amino acids $V$, the mask $\mathsf m$ and the
special tokens) and $f_\phi$ the logits of the frozen model \citep{lin2023esm}.
Fix a state $x$, a masked position $i$, and write $s=x_i=\mathsf m$ for the
current token. With the inverse temperature $\beta_t$ of the reference
($\beta_t=1-t$ in our implementations of both arms, where $t$ is the masking
level and $t=1$ the fully masked sequence),
\begin{equation}\label{eq:token-kernel}
r(y)=\operatorname{softmax}_y\!\big(\beta_t f_\phi(x)_i\big),\qquad
p_\theta(y)=\operatorname{softmax}_y\!\big(\beta_t f_\phi(x)_i+u_\theta(x,y,t)\big),\qquad y\in\mathcal T,
\end{equation}
where $u_\theta(x,\cdot,t)\in\mathbb R^{\mathcal T}$ is the control at the edges
leaving $x$ at position $i$. Both are probability vectors, so the token kernel is
the rate~\eqref{eq:ref-rates} at unit intensity; $\lambda_t$ enters only through
the sampler. Let $q(a)=\mathbb P(X^i_1=a\mid X_t=x)$ be the posterior of the clean
token, supported on $V$. Training minimises the cross-entropy
\begin{equation}\label{eq:ce}
\mathrm{CE}(u)=-\sum_{y\in \mathcal T} q(y)\log p_\theta(y)
\end{equation}
(as in MadSBM, the loss is also evaluated at unmasked positions, where it is a
copy task; the analysis concerns masked positions, the only ones where the
sampler uses the kernel). Its gradient in the logits is
$p_\theta(y)-q(y)$, which sums to zero: the gradient has no component along
$\mathbf 1$, so nothing in the cross-entropy ever selects the gauge.

\begin{definition}[Vertex-potential control]\label{def:vertex}
For a potential $V:\mathcal X\to\mathbb R$ the discrete gradient is the edge
function $(\bar\nabla V)(x,y)=V(y)-V(x)$, i.e.\ the difference of the potential
between the two ends of an edge. Since $V$ cannot be evaluated on all of
$\mathcal X$, we read it locally: at a state $x$, a time $t$ and a masked
position $i$, an energy head outputs
$\tilde v_\theta(x,t)\in\mathbb R^{\mathcal T}$ whose coordinate $a$ is the
potential of the neighbour $x^{i\leftarrow a}$. As $x^{i\leftarrow x_i}=x$, the
coordinate $x_i$ is the potential of the current state itself, so
\begin{equation}\label{eq:vertex-control}
u_\theta\big(x,\,x^{i\leftarrow a},\,t\big)=\tilde v_\theta(x,t)_a-\tilde v_\theta(x,t)_{x_i} .
\end{equation}
The head is initialised at zero, so $u_\theta\equiv0$ and $p_\theta=r$ at
initialisation: the untrained model is the tempered protein language model.
\end{definition}

By construction $u_\theta(x,s,t)=0$: the current token carries zero control. We
call this the \emph{vertex gauge}. The implemented action penalty is the
per-state sum
\begin{equation}\label{eq:apsi}
\mathcal A_\Psi(u;x)=\sum_{y\in\mathcal T}r(y)\Psi\big(u(y)\big),\qquad \Psi(z)=e^z-z-1,
\end{equation}
which matches the kernel of MadSBM's action metric at unit intensity $\nu=1$.
MadSBM trains with the cross-entropy alone; GDSBM adds the action to the loss and
minimises
\begin{equation}\label{eq:gdsbm-loss}
\mathcal L(\theta)=\mathbb E\big[\mathrm{CE}\big]+\lambda\,\mathbb E\big[\mathcal A_\Psi\big],
\end{equation}
with $\mathcal A_\Psi$ evaluated on the vertex-gauge
control~\eqref{eq:vertex-control} and $\lambda\ge0$ ($\lambda=0.1$ for the GDSBM
arm in all experiments).

\paragraph{Notation.} Two control classes act on the logits at a masked position
with current token $s$: the free class
$\mathcal U_{\rm free}=\mathbb R^{\mathcal T}$ (MadSBM) and the vertex class
$\mathcal U_{\rm vtx}=\{u\in\mathbb R^{\mathcal T}:u_s=0\}$. A control $u$ induces
$p_\theta[u](a)=r(a)e^{u_a}/Z_u$ with $Z_u=\sum_{b}r(b)e^{u_b}$. Let $q$ be the
clean-token posterior, supported on $V$, and set $\nu_V=\sum_{a\in V}r(a)$,
$\tilde r=r|_V/\nu_V$ and $J(p)=\mathrm{CE}(p)+\lambda\mathcal A_\Psi$, with
$\mathcal A_\Psi$ read through~(b) below as a function of $p$. The result states
three facts: the vertex class loses no conditional (a); in that class the action
is a function of the conditional alone (b); and the penalty has an exact effect
on the minimiser (c).

\begin{theorem}[Vertex gauge: expressiveness, anchor, shrinkage]\label{thm:vertex}
\begin{enumerate}[label=(\alph*),leftmargin=1.8em,itemsep=0.1em,topsep=0.2em]
\item \emph{(pinning the gauge costs nothing)} $p_\theta[u]=p_\theta[w]$ if and
only if $u-w\in\mathbb R\mathbf 1$. Hence every orbit meets
$\mathcal U_{\rm vtx}$ exactly once, at $u-u_s\mathbf 1$. Moreover
$u\in\mathcal U_{\rm vtx}$ if and only if $u_a=v_a-v_s$ for some
$v\in\mathbb R^{\mathcal T}$, and a conditional $p$ is recovered by
$u_a=\log\big[(p_a/r(a))/(p_s/r(s))\big]$. Consequently the two classes realise
the same conditionals and the same infimum of~\eqref{eq:ce}, so GDSBM at
$\lambda=0$ is functionally identical to MadSBM. What the conditionals do not fix
is the rate level: $R^0e^{u}$ and $R^0e^{u-u_s\mathbf 1}$ differ by the factor
$e^{-u_s}$ isolated in Theorem~\ref{thm:gauge};
\item \emph{(the action becomes a model quantity)} for $u\in\mathcal U_{\rm vtx}$ with $p=p_\theta[u]$,
\begin{equation}\label{eq:vgauge}
\mathcal A_\Psi(u;x)=\operatorname{KL}(r\,\|\,p)+h\!\Big(\frac{r(s)}{p(s)}\Big),
\end{equation}
a function of $p$ alone, whose second term vanishes if and only if $p(s)=r(s)$;
\item \emph{(exact effect of the penalty)} assume $r(y)>0$ for all $y$ and let
the class be expressive, i.e.\ any interior $p$ is attainable (by~(a) this holds
in either class). Then $J$ has an interior minimiser, and every minimiser
satisfies, on the amino acids,
\begin{equation}\label{eq:shrink}
p(\cdot\mid V)=(1-w)\,q+w\,\tilde r,\qquad w=\frac{\lambda\nu_V}{1+\lambda\nu_V}.
\end{equation}
In particular $\|p(\cdot\mid V)-q\|_{\rm TV}=w\,\|\tilde r-q\|_{\rm TV}$.
\end{enumerate}
\end{theorem}

The proof is in Appendix~\ref{app:vertex}. Part~(b) is
Theorem~\ref{thm:gauge}(a) at $\nu=1$ with $\zeta$ eliminated through
$p(s)=r(s)/\zeta$; part~(c) is a stationarity computation on the open simplex, in
which the penalty turns the loss into a cross-entropy with weights
$q+\lambda r$. State by state, the penalised minimiser is thus a \emph{convex
shrinkage} of the data posterior towards the amino-acid reference, with weight
$w\simeq8.5\%$ at $\lambda=0.1$ and $\nu_V\simeq0.93$: an explicit, measurable
regulariser towards the ESM-2 prior, and nothing more
(Remarks~\ref{rem:shrink}--\ref{rem:geometry}).

\begin{table}[t]
\centering\small
\caption{\textbf{Form $\times$ penalty}, seed $42$, $300$ steps, common random numbers. \emph{val CE} is the cross-entropy on the validation split; \emph{action} is the training-convention action (raw gauge, averaged over tokens) of the control as parameterised.}
\label{tab:ablation}
\begin{tabular}{@{}llcc@{}}
\toprule
Class & $\lambda$ & val CE & action\\
\midrule
free (MadSBM) & $0$ & $1.36985924$ & $15.73$\\
vertex (ours) & $0$ & $1.36985921$ & $42.91$\\
\midrule
free (MadSBM) & $0.1$ & $1.4510$ & $0.413$\\
vertex (ours) & $0.1$ & $1.4473$ & $0.584$\\
\bottomrule
\end{tabular}
\end{table}

\paragraph{Training.} At each step a masking level $t$ is drawn on a uniform grid
and each residue is masked independently with probability $t$; the frozen
language model gives $f_\phi(x_t)$, the head gives $\tilde v_\theta$, and the
loss~\eqref{eq:gdsbm-loss} is computed with $\mathcal A_\Psi$ on
$u_\theta=\bar\nabla\tilde v_\theta$. Only the head is trained (hyperparameters in
Appendix~\ref{app:details}). In Algorithm~\ref{alg:train}, $\mathcal P$ is the set
of residue positions of the batch: the loss denoises at masked positions and
copies at unmasked ones. The two arms differ only at line~\ref{line:arm}, where
the head output is read as a vertex potential (GDSBM) or as the control itself
(MadSBM, with $\lambda=0$).

\paragraph{Sampling.} We use MadSBM's sampler unchanged, with the same
engineering constants, for both models: only the control field differs, so the
comparison bears on a single object, the control class and the penalty.
Algorithm~\ref{alg:sample} starts from the fully masked sequence and, at each of
$N$ steps, uses the control to decide which masked positions jump and which amino
acid they receive. Two features bear on what the audit of
Section~\ref{sec:ex-audit} measures. The jump probability uses the control damped
by a rate scale $\kappa=0.01$, whereas the token kernel uses it undamped, so the
two uses of $u_\theta$ are not on the same scale; in particular the jump schedule,
unlike the token kernel, is not gauge-invariant. And the candidate set is
restricted to the $20$ amino acids after the nucleus filter, so the tokens
carrying the mass $1-\nu_V$ of Theorem~\ref{thm:vertex}(c) are never emitted.
Parameters $\kappa,\eta,\tau,\pi_{\rm top}$ are in Appendix~\ref{app:details}.
Guidance (Section~\ref{sec:conditioning}) adds $\beta\Delta\Phi$ to $u$ before
lines~\ref{line:intensity} and~\ref{line:kernel}, and nowhere else.

\begin{algorithm}[t]
\caption{GDSBM training step (MadSBM: free reading, $\lambda=0$)}
\label{alg:train}
\begin{algorithmic}[1]
\Require corpus $\mathcal D$, frozen logits $f_\phi$, head $\theta$, penalty weight $\lambda$, grid $\mathcal G$
\State draw $x_1\sim\mathcal D$, $t\sim\mathrm{Unif}(\mathcal G)$; mask each residue with probability $t$ to get $x_t$
\State $r_i\gets\operatorname{softmax}\big(\beta_tf_\phi(x_t)_i\big)$, \quad $\beta_t=1-t$
\State $\tilde v\gets\mathrm{head}_\theta(x_t,t)$; \quad $u_i(a)\gets\tilde v_i(a)-\tilde v_i\big((x_t)_i\big)$ \;\emph{(GDSBM)}, \quad $u_i(a)\gets\tilde v_i(a)$ \;\emph{(MadSBM)} \label{line:arm}
\State $p_i\gets\operatorname{softmax}\big(\beta_tf_\phi(x_t)_i+u_i\big)$
\State $\theta\gets\mathrm{AdamW}\big(\theta,\ \nabla_\theta[\mathrm{CE}(p,x_1)+\lambda\,\mathcal A_\Psi(u;r)]\big)$
\end{algorithmic}
\end{algorithm}
Algorithm~\ref{alg:sample} in Appendix~\ref{app:details} records the sampler.

\paragraph{Scope.} The word \emph{bridge} is used in the sense of the Background:
a jump process on the sequence graph, anchored at a fully masked source, with a
reference given by the tempered protein language model and a control of the
form~\eqref{eq:doob-control}. GDSBM is a parametric approximation of this object,
with locally read potentials (Remark~\ref{rem:geometry}), trained by a
regularised per-state objective. We claim exactly
Theorems~\ref{thm:gauge}--\ref{thm:cond}, and do not claim that the trained model
minimises the path-space entropy of the Schr\"odinger problem.

\section{Conditioning: Algebra on Potentials}\label{sec:conditioning}

Design asks for sequences with prescribed properties. The vertex form is the
level at which conditioning composes, and on peptides it is also the level at
which the bias of practical guidance becomes measurable rather than asserted.

\paragraph{Three objects.} Let $\mathbb P^\ast$ be the bridge of
Section~\ref{sec:background}, with terminal law $\mu_1$, and let
$\Phi:\mathcal X\to\mathbb R$ be a bounded design objective. The \emph{tempered
target} at inverse temperature $\beta\ge0$ is
\begin{equation}\label{eq:tempered}
\mu_1^\beta(x)=\frac{\mu_1(x)\,e^{\beta\Phi(x)}}{Z_\beta},\qquad Z_\beta=\sum_x\mu_1(x)e^{\beta\Phi(x)} .
\end{equation}
Computing the tilt potential
$v^{\beta\Phi}_t(x)=\log\mathbb E_{\mathbb P^\ast}[e^{\beta\Phi(X_T)}\mid X_t=x]$
requires an expectation over continuations, intractable on $\mathcal X$; the
standard surrogate replaces it by the reward $\beta\Phi$ itself. We evaluate
$\Phi$ on partially masked sequences by letting masked positions contribute
nothing and write, for a state $x$, a position $i$ and a token $a\in V$,
$\Delta\Phi(x,i,a)=\Phi(x^{i\leftarrow a})-\Phi(x)$. The \emph{myopic guided
control} is $u^\beta_\theta=\bar\nabla(v_\theta+\beta\Phi)$, that is $u_\theta$
plus $\beta\Delta\Phi$ on every edge, and the token
kernel~\eqref{eq:token-kernel} becomes
\begin{equation}\label{eq:myopic-kernel}
p^\beta_\theta(a)\;\propto\;p_\theta(a)\,e^{\beta\Delta\Phi(x,i,a)},\qquad a\in V .
\end{equation}
The same increment enters the jump intensity and the token logits, so the guided
sampler is one modification of the control rather than two; this is the discrete
counterpart of classifier guidance
\citep{dhariwal2021diffusion,ho2022classifier,nisonoff2024,schiff2024,rectorbrooks2024}.
Finally, let $x^{(0)}\to\dots\to x^{(N)}$ be the sampler's states, $x^{(0)}$ the
fully masked source and $x^{(N)}$ the fully decoded sequence, and run $K$
particles with the \emph{unguided} kernel as proposal and incremental weights
\begin{equation}\label{eq:incr-weight}
G_k\big(x^{(k-1)},x^{(k)}\big)=\exp\!\Big(\beta\big[\Phi(x^{(k)})-\Phi(x^{(k-1)})\big]\Big),\qquad k=1,\dots,N ,
\end{equation}
with self-normalised weights $W^j\propto\prod_{k=1}^N G^{\,j}_k$. For a statistic
$c:\mathcal X\to\mathbb R$ with finitely many values and
$\psi:\mathbb R\to\mathbb R$, write $\Phi=\psi\circ c$. The result states three
facts: (a)~a tilt is an $h$-transform whose potential \emph{adds} to the trained
one, with an explicit gap when two tilts are composed; (b)~the tempered target
has a closed-form marginal for any scalar statistic; (c)~the myopic sampler is
biased, while the weighted particle system recovers the target.

\begin{theorem}[Conditioning]\label{thm:cond}
\begin{enumerate}[label=(\alph*),leftmargin=1.8em,itemsep=0.1em,topsep=0.2em]
\item \emph{(algebra on potentials)} The path law $\mathbb P^{\beta}$ with rates
$R^{\ast}_t(x,y)\,e^{v^{\beta\Phi}_t(y)-v^{\beta\Phi}_t(x)}$ has terminal law
$\mu_1^\beta$, and its control relative to the reference is the vertex potential
\begin{equation}\label{eq:add-control}
\begin{split}u^{\beta}_t(x,y)={}&\big(v^{\ast}_t(y)-v^{\ast}_t(x)\big)\\&+\big(v^{\beta\Phi}_t(y)-v^{\beta\Phi}_t(x)\big).\end{split}
\end{equation}
For bounded $\Phi_1,\Phi_2$, tilting by $\Phi_1$ then $\Phi_2$ gives
$\mathbb P^{\Phi_1+\Phi_2}$, but the potentials add only up to a gap,
$v^{\Phi_1+\Phi_2}=v^{\Phi_1}+v^{\Phi_2}+G$, with $G(\cdot,T)=0$ and
\begin{equation}\label{eq:compose-gap}
\begin{split}
G(x,t)={}&\log\mathbb E_{\mathbb P^{\Phi_1}}\big[e^{\Phi_2(X_T)}\mid X_t=x\big]\\
&-\log\mathbb E_{\mathbb P^{\ast}}\big[e^{\Phi_2(X_T)}\mid X_t=x\big].
\end{split}
\end{equation}
\item \emph{(closed-form target)} For every value $v$ of $c$, with no product
structure on $\mu_1$, no fixed length and no additivity of $c$ assumed,
\begin{equation}\label{eq:exact-marginal}
\mu_1^\beta(c=v)=\mu_1(c=v)\,e^{\beta\psi(v)}/Z_\beta,\qquad Z_\beta=\textstyle\sum_{v'}\mu_1(c=v')\,e^{\beta\psi(v')}.
\end{equation}
\item \emph{(bias and its removal)} The myopic kernel~\eqref{eq:myopic-kernel} is
exact for the one-step tilt and biased for the terminal tilt~\eqref{eq:tempered};
the weights telescope, $\prod_{k}G_k=e^{\beta[\Phi(x^{(N)})-\Phi(x^{(0)})]}$ with
$\Phi(x^{(0)})$ constant, so $W^j\propto e^{\beta\Phi(x^{(N),j})}$ and
$\sum_j W^j\delta_{x^{(N),j}}$ is a consistent estimator of $\mu_1^\beta$ as
$K\to\infty$, with or without intermediate resampling.
\end{enumerate}
\end{theorem}

The proof is in Appendix~\ref{app:cond}: (a) is Fact~\ref{fact:gibbs} with
$\mathbb P^\ast$ in place of $\mathbb P^0$, plus a change of measure for the gap;
(b) groups~\eqref{eq:tempered} by the level sets of $c$; (c) is a telescoping
product. Since a sum of vertex potentials is one, the conditioned model stays in
the class of Definition~\ref{def:vertex}, whereas on a free edge field the tilt
inherits an unidentified shift. The gap vanishes only at $t=T$, so the additive
rule used in practice is the \emph{myopic} approximation, not an identity; what
(a) does give at every $\beta$ is that the \emph{target} of a multi-objective
sweep is a single tilted law, so a sweep needs no retraining
(Remark~\ref{rem:cond-scope} delimits what we measure).

\paragraph{The bias becomes a measurement.} Peptide design conditions are scalar
statistics of the sequence: the net charge at pH~7,
\begin{equation}\label{eq:charge}
\begin{gathered}
c(x)=\sum_{i=1}^{L}\varphi(x_i),\\
\varphi(\mathrm K)=\varphi(\mathrm R)=+1,\ \ \varphi(\mathrm D)=\varphi(\mathrm E)=-1,\ \ \varphi(\mathrm H)=+0.1,
\end{gathered}
\end{equation}
and $\varphi=0$ otherwise, or the mean hydropathy over decoded residues
\citep{kyte1982}, which is not additive in the positions. A target is a function
of such a statistic, $\Phi=\psi\circ c$; we use $\psi(v)=-(v-c^\star)^2$. Part~(b)
then gives the one-dimensional marginal of the tempered law in closed form, so
the bias of the myopic sampler documented in
\citep{lee2025debiasing,wan2026dgm} can be compared to its target rather than
asserted.

\begin{remark}[The measurement protocol]\label{rem:protocol}
Throughout, $\mu_1$ is the terminal law of the sampler actually run at $\beta=0$,
with its temperature and nucleus truncation, not an idealised bridge and not the
data. We estimate its marginal $p(v)=\mu_1(c=v)$ from $n_0$ unguided samples,
$\hat p$, and form $q_\beta(v)\propto\hat p(v)\,e^{\beta\psi(v)}$ on the support
of $\hat p$. For $n$ guided samples with empirical charge histogram
$\hat\mu_\beta$ we report
$\mathrm{TV}(\hat\mu_\beta,q_\beta)=\tfrac12\sum_{v}|\hat\mu_\beta(v)-q_\beta(v)|$,
the sum running over the union of the supports, together with the fraction of
samples in a window $[c^\star-1,c^\star+1]$ against its value under $q_\beta$. A
value of $\mathrm{TV}$ is read against a \emph{sampling floor}: its law under the
null in which the sampler draws exactly from the true tempered marginal, obtained
by a parametric bootstrap from $\hat p$ that redraws both $\hat p$ ($n_0$
samples) and the guided sample ($n$ samples). Two limits apply: matching the
marginal of $c$ is necessary for $\hat\mu_\beta$ to approximate $\mu_1^\beta$ but
not sufficient, since the statistic sees one coordinate of the sequence only; and
$\hat p$ is itself an estimate, whose error the floor accounts for.
\end{remark}

\section{Experiments}\label{sec:experiments}

The experiments answer three questions, in order. Is the action of a
cross-entropy-trained control a statement about transport, or about the gauge
(Section~\ref{sec:ex-audit})? Does the vertex form change the model when the
penalty is off, and what does the penalty do when the form is free
(Section~\ref{sec:ex-ablation})? And does the regularised model pay for it in
sequence quality (Section~\ref{sec:ex-uncond})? Section~\ref{sec:ex-cond} reports
the conditioning measurement. Deviations and threats to validity are listed in
Appendix~\ref{app:details}, and full tables in Appendix~\ref{app:tables}.

\paragraph{Setup.} \emph{MadSBM (reimpl.)} is the free edge field trained by
cross-entropy alone ($\lambda=0$); \emph{GDSBM (ours)} is the vertex-potential
control of Definition~\ref{def:vertex} trained with~\eqref{eq:gdsbm-loss} at
$\lambda=0.1$. They share the frozen ESM-2 650M backbone \citep{lin2023esm}, the
AdaLN head ($\approx45$M trainable parameters), corpus, optimiser, schedule
($127\,000$ steps) and sampler, so they differ in exactly two places: the control
class and the penalty. The dataset of \citep{goel2026madsbm} is not public, so we
rebuild a corpus from UniProt (reviewed, $5\le L\le50$, strict 20-letter
alphabet) split $90/5/5$ at seed $42$: all comparisons are internal to our two
arms. We train three seeds per arm ($42,43,44$) and sample at
$N\in\{32,64,128\}$, generating $20$ sequences at each length
$L\in\{5,\dots,50\}$ per (arm, seed, $N$), i.e.\ $920$ sequences and $16\,560$ in
all; pairs are always at equal seed and $N$. We report pseudo-perplexity under
ESM-2, mean pLDDT from ESMFold rescaled to $[0,100]$, the fraction of unique
sequences (Appendix~\ref{app:details}), and the action in three conventions---raw,
vertex gauge, and gauge-fixed, the minimum over the gauge class, equal to
$\operatorname{KL}(r\|p)$ by Theorem~\ref{thm:gauge}---plus its forward-kernel
counterpart.

\subsection{The raw action of a free control measures its gauge}\label{sec:ex-audit}

\begin{figure}[t]
\centering
\includegraphics[width=0.95\linewidth]{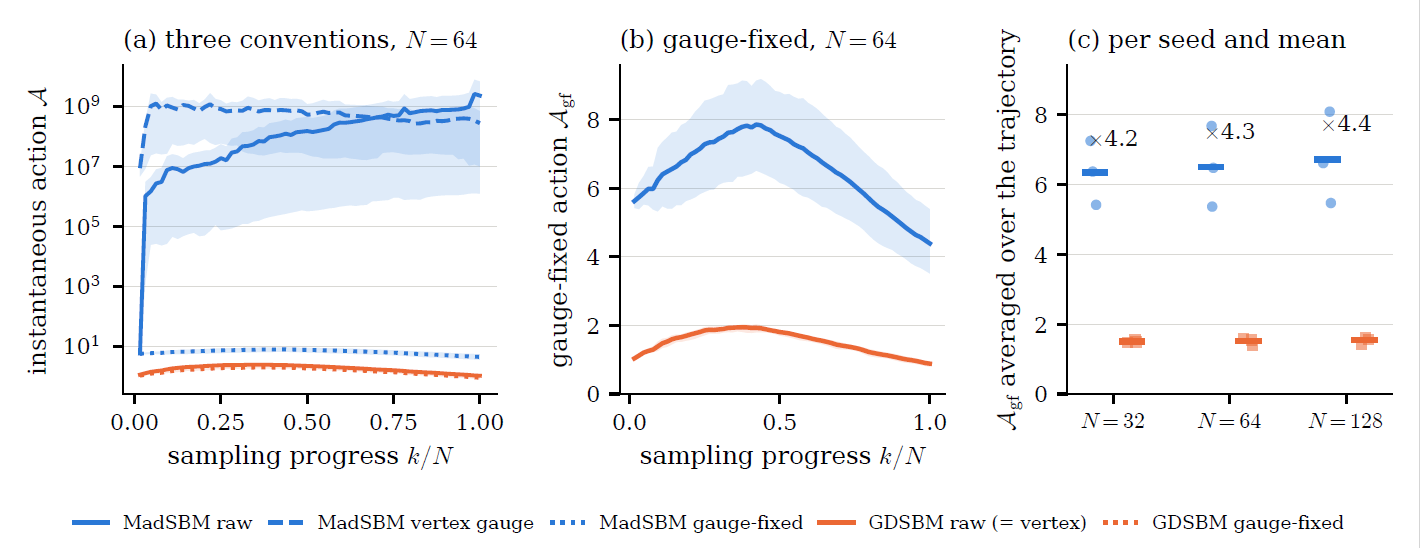}
\caption{\textbf{The raw action of a free control measures its gauge.} Seed-averaged action along sampling trajectories ($33$ tokens, three seeds, $46$ sequences per arm, seed and $N$). \emph{(a)}~Three conventions at $N\in\{32,64,128\}$ (log scale): the raw and vertex-gauge actions of MadSBM are of order $10^{8}$, those of GDSBM stay near $1.8$. \emph{(b)}~After gauge fixing, $\mathcal A_{\rm gf}=\operatorname{KL}(r\|p)$ is a property of the conditional alone: GDSBM is lower by a factor $4.2$--$4.4$ (error bars: inter-seed s.d.). Source: Table~\ref{tab:gauge}.}
\label{fig:action-gauge}
\end{figure}

Figure~\ref{fig:action-gauge} evaluates the instantaneous action along the
sampling trajectories; numerical values are in Table~\ref{tab:gauge}. Three
readings, in increasing order of what they license.

\emph{(i)~The raw action of the free control is not a transport quantity.} For
MadSBM it is of order $10^{8}$, with a seed-to-seed standard deviation larger
than its mean and a spread of a factor $3.6\cdot10^{3}$ across the three seeds
(from the logs), while the corresponding gauge-fixed value is $6.35$--$6.72$
across the three budgets, with a seed-to-seed standard deviation of
$14$--$20\%$ of its mean. Theorem~\ref{thm:gauge} identifies the discrepancy: the raw
action is $\operatorname{KL}(r\|p)+h(\zeta)$ and the cross-entropy never
constrains $\zeta$, so the second term is free to be large. This is a statement
about the quantity, not about the model: two controls differing by a gauge induce
the same token kernel at every state. We therefore regard any comparison of raw
actions between a free field and a pinned one as void, including the comparison
in our own favour.

\emph{(ii)~After gauge fixing, the comparison is well posed and GDSBM has the
lower action.} The gauge-fixed action is a function of the learned conditional
alone. GDSBM is below MadSBM by a factor $4.2$, $4.3$ and $4.4$ at
$N=32,64,128$ on the full vocabulary, and $3.8$, $3.9$ and $4.0$ on the amino
acids, with the ordering holding for all nine (seed, $N$) pairs. As a numerical
check, the gauge-fixed column of Table~\ref{tab:gauge} reproduces the
independently computed $\operatorname{KL}(r\|p)$ to the printed precision, as
Theorem~\ref{thm:gauge}(a) requires.

\emph{(iii)~At the forward kernel, the gap narrows sharply.}
Figure~\ref{fig:action-gauge} uses the reverse kernel $\Psi$, which is what the
penalty applies. Auditing the same checkpoints at the Girsanov kernel $a$, with
the closed-form gauge minimisation of Theorem~\ref{thm:gauge}(c), the raw forward
action of MadSBM lies between $7\cdot10^{9}$ and $2.4\cdot10^{10}$, but the
\emph{gauge-fixed} forward action of the two arms differs by a factor of only
$1.07$--$1.09$, with the sign holding for all nine (seed, $N$) pairs. This is not
saturation: both values are around $0.51$--$0.56$, well below the bound $1$ of
Theorem~\ref{thm:gauge}(c) at unit intensity. Inverting its closed form,
$\operatorname{KL}(p\|r)=-\log(1-\mathcal A^{a}_{\rm gf})$, on the seed means
gives a ratio of $1.11$--$1.13$ between the arms, against $4.2$--$4.4$ for
$\operatorname{KL}(r\|p)$. The reduction is thus concentrated in the reverse
direction, the one the penalty acts on. We report both (Table~\ref{tab:fwd}) and
make no claim that it transfers quantitatively to the path-space action.

\subsection{Form and penalty, separated}\label{sec:ex-ablation}

Table~\ref{tab:ablation} switches the two changes independently, at seed $42$ for
$300$ steps with identical data order and noise draws. \emph{At $\lambda=0$ the
two forms are the same model, to eight digits:} the validation cross-entropies
agree to $3.6\cdot10^{-8}$, i.e.\ to float32 non-determinism, which is
Theorem~\ref{thm:vertex}(a) measured. Yet the raw actions of that \emph{same}
model differ by a factor $2.7$, purely because the two parameterisations pick
different representatives of the gauge class; a quantity that moves by $2.7$
while the model is unchanged to eight digits cannot be evidence of minimal-action
transport. \emph{At $\lambda=0.1$ the penalty is partly absorbed by the gauge of
a free field:} the free field reaches a \emph{lower} raw action ($0.413$ against
$0.584$) at a \emph{worse} cross-entropy ($1.4510$ against $1.4473$). By
Theorem~\ref{thm:gauge}(b) it can drive $h(\zeta)$ to zero at no cost in
cross-entropy, so part of the penalty's gradient is spent on the gauge, whereas
the vertex form must pay the anchor $h(r(s)/p(s))$ of
Theorem~\ref{thm:vertex}(b), reducible only by moving $p$. The penalty therefore
does not mean the same thing in the two classes, which is why we treat the form
and the penalty as a single design choice. These runs are $300$ steps at one
seed: a mechanism check, not a performance comparison.

\subsection{Unconditional generation at matched compute}\label{sec:ex-uncond}

Across nine paired tests (three metrics, three budgets, three seeds, $2$ degrees
of freedom each; Table~\ref{tab:uncond}), no difference survives a Bonferroni
threshold of $0.05/9=0.0056$. The smallest $p$-value is pLDDT at $N=64$,
$\Delta=-0.24$ with $p=0.014$ uncorrected and a constant sign across the three
seeds; perplexity changes sign with the budget ($+0.25$, $-0.04$, $+0.04$);
uniqueness favours GDSBM with a constant sign at $N=32$ and $64$ ($+0.008$,
$p=0.073$; $+0.010$, $p=0.118$) and not at $N=128$. With three seeds a $p$-value
above $0.05$ is not evidence of equivalence, and the sign pattern is the
informative part. The conclusion we draw is the negative one we need: the
four-fold reduction of the gauge-fixed action of Section~\ref{sec:ex-audit} is
\emph{not} bought at a detectable cost in sequence quality, with a possible small
pLDDT penalty of about $0.1$--$0.5$ points that these three seeds can neither
establish nor exclude. Checkpoint cross-entropies, re-evaluated with common
random numbers, are in Table~\ref{tab:stability} and show a paired difference of
alternating sign.

\subsection{Conditioning against a closed-form target}\label{sec:ex-cond}

We guide each sampler with the myopic twist of~\eqref{eq:myopic-kernel}, using
the reward $\Phi=-(c-c^\ast)^2$ of~\eqref{eq:charge} on the net charge $c$, at
tempering $\beta\in\{0.5,1,2,4\}$ with $c^\ast=3$. By Theorem~\ref{thm:cond} the
tempered target law of $c$ is available in closed form,
$q_\beta(v)\propto \hat p(v)\,e^{-\beta (v-c^\ast)^2}$, where $\hat p$ is the
unguided law of $c$ estimated from $n_0=300$ unconditional samples. We report
$\mathrm{TV}(\hat\mu_\beta,q_\beta)$ between the empirical law of $c$ under
guidance ($n=256$ samples) and $q_\beta$, and the empirical mass of the window
$c\in[2,4]$ against the mass $q_\beta$ assigns to it; each arm is compared with
its \emph{own} target, since $\hat p$ differs between arms. Even a perfect
sampler has $\mathrm{TV}>0$ because $\hat p$ and $\hat\mu_\beta$ are
finite-sample estimates; we quantify this \emph{sampling floor} by a parametric
bootstrap (5000 replicates, same $n_0$ and $n$): draw
$\hat p^{(r)}\sim\mathrm{Mult}(n_0,\hat p)/n_0$, form its tilt, and compare with
an independent $\mathrm{Mult}(n,q_\beta)/n$ draw. The floor is the mean of the
resulting TV, averaged over seeds; mass of $\hat p$ outside its support is
ignored, so it is slightly underestimated.

\begin{table}[t]
\centering\footnotesize
\caption{Conditioning on net charge, three seeds ($42,43,44$; one checkpoint per seed and arm). Entries are mean$\pm$inter-seed s.d., followed after the slash by the sampling floor (TV) or the target mass $q_\beta([2,4])$ (window). Window rates in \%. Bold marks the better arm at each $\beta$; no significance is claimed (see text).}
\label{tab:cond}
\resizebox{\textwidth}{!}{%
\begin{tabular}{@{}lcccc@{}}
\toprule
 & \multicolumn{2}{c}{TV to $q_\beta$ / floor $\downarrow$} & \multicolumn{2}{c}{window rate / target (\%) $\uparrow$}\\
\cmidrule(lr){2-3}\cmidrule(lr){4-5}
$\beta$ & MadSBM & GDSBM & MadSBM & GDSBM\\
\midrule
0.5 & 0.549$\pm$0.098 / 0.081 & \textbf{0.473$\pm$0.023} / 0.087 & 56.0 / 88.5 & \textbf{68.6} / 86.8 \\
1   & 0.470$\pm$0.064 / 0.075 & \textbf{0.380$\pm$0.057} / 0.080 & 73.7 / 96.9 & \textbf{81.0} / 96.7 \\
2   & 0.350$\pm$0.052 / 0.061 & \textbf{0.253$\pm$0.058} / 0.066 & 88.4 / 99.4 & \textbf{93.0} / 99.5 \\
4   & 0.281$\pm$0.054 / 0.027 & \textbf{0.178$\pm$0.077} / 0.044 & 97.0 / 99.9 & \textbf{99.5} / 99.9 \\
\bottomrule
\end{tabular}}
\end{table}

\begin{figure}[t]
\centering
\includegraphics[width=0.95\linewidth]{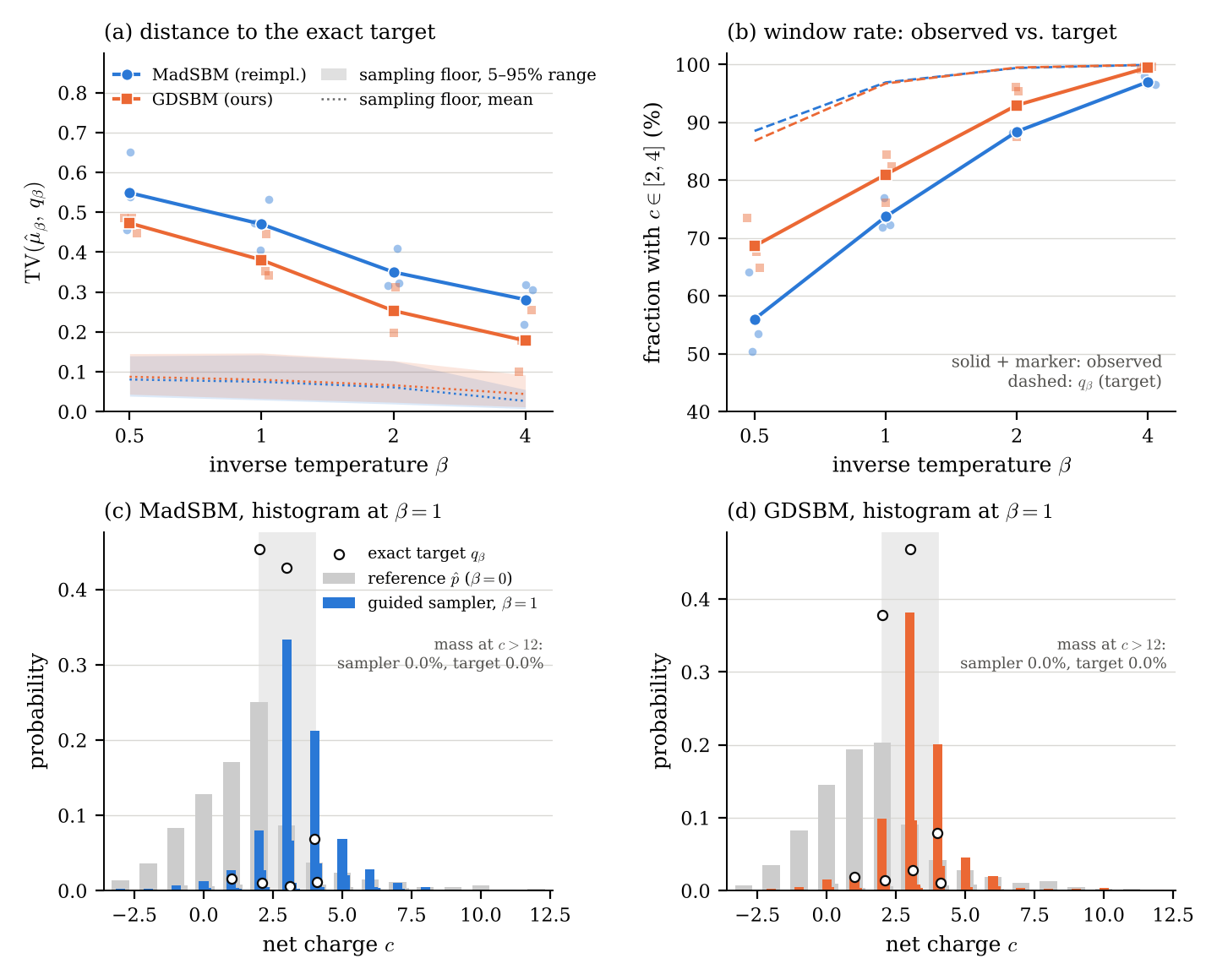}
\caption{Conditioning on net charge, three seeds ($42,43,44$). \textbf{(a)}~$\mathrm{TV}(\hat\mu_\beta,q_\beta)$ versus $\beta$: large markers are seed means, small markers individual seeds; the band is the 5--95\% range of the sampling floor, the dotted line its mean. \textbf{(b)}~Empirical mass of the window $c\in[2,4]$ (solid) against the target $q_\beta([2,4])$ (dashed). \textbf{(c,d)}~Law of $c$ at $\beta=1$, averaged over seeds (bars), against the unguided $\hat p$ (grey) and the closed-form target $q_\beta$ (circles), for MadSBM~(c) and GDSBM~(d); the strip is the window.}
\label{fig:conditioning}
\end{figure}

\emph{(i)~Both samplers move toward the target, neither reaches it.} The observed
TV exceeds the 95\% quantile of its own sampling floor in all $24$ (arm, $\beta$,
seed) cases, one marginally (GDSBM, seed~43, $\beta=4$: $0.100$ against $0.099$),
and the mean TV exceeds the mean floor by a factor of at least $3.8$ (smallest:
$0.253/0.066$, GDSBM at $\beta=2$). \emph{(ii)~The window rate is below its target in all $24$
cases}, the expected gap of a twist that sees only one-step reward increments,
not a property of either model. \emph{(iii)~Ordering.} GDSBM has the lower TV in
$11$ of the $12$ paired (seed, $\beta$) comparisons and the higher window rate in
$11$ of $12$ (exceptions: seed~43 at $\beta=0.5$, seed~44 at $\beta=2$), with
mean TV gaps $0.076$, $0.090$, $0.097$, $0.103$, the first within MadSBM's
inter-seed spread ($0.098$). With three seeds we report direction, not
significance, and do not attribute the ordering to the vertex form rather than
the penalty, the $2\times2$ ablation not having been run under guidance. One
checkpoint per seed and arm, evaluated on CPU; the particle system of
Theorem~\ref{thm:cond}(c) is implemented but not evaluated here.

\section{Conclusion}\label{sec:conclusion}

A discrete bridge trained by cross-entropy leaves one degree of freedom per state
unconstrained, and the action it reports lives almost entirely there.
Theorem~\ref{thm:gauge} makes this exact and gives the comparable, gauge-fixed
quantity; Theorem~\ref{thm:vertex} pins the gauge by the form of the Doob
optimum, at no cost in representable conditionals, and turns an action penalty
into a convex shrinkage towards the reference; Theorem~\ref{thm:cond} makes
conditioning addition of potentials and replaces an asserted bias by a measured
one. Two questions remain open (Appendix~\ref{app:future}): a per-state
regulariser is not yet a path-space action, and a locally read potential is not
yet a global one. We plan to study the geometry of the control potential through
the Onsager-type and entropy-dissipation structure of finite-state
Schr\"odinger bridges \citep{benabdallah2026nonlocal}.

\appendix

\section{Deferred remarks}\label{app:remarks}

This appendix collects the remarks referenced from the main text.

\begin{remark}[The reported action is not a property of the model]\label{rem:gauge}
Cross-entropy identifies a control only up to the gauge, because the softmax
cannot see it, yet the action is not gauge-invariant. The number a model reports
as its action is therefore the sum of a genuine property of the model, the
divergence $\nu\operatorname{KL}(r\|p)$, and a term that by
Theorem~\ref{thm:gauge}(b) can be made arbitrarily large without altering the
token kernel at any state. Any comparison of raw actions is thus void, including
one in our own favour: two parameterisations that agree as models to eight digits
report raw actions differing by a factor of $2.7$, and a trained free field
reports $10^{8}$ for a gauge-fixed value of $6.35$--$6.72$
(Section~\ref{sec:ex-audit}). The two directions are moreover not symmetric:
$h(s)\sim s$ for large $s$ and $h(s)\sim\log(1/s)$ for small $s$, so inflating the
action is exponentially cheap in the shift while deflating it is only
logarithmic. The vertex form of Section~\ref{sec:gdsbm} does not remove the second
term; it pins it, by imposing $u(x_i)=0$ at the current token. That term becomes
the anchor $h(r(x_i)/p(x_i))$, a function of the learned conditional alone, so the
action is again a property of the model and a penalty acts on the conditional
rather than on the gauge.
\end{remark}

\begin{remark}[Two kernels, two scales]\label{rem:forward-gauge}
Theorem~\ref{thm:gauge}(c) is the same computation at the forward kernel $a$ of
Fact~\ref{fact:girsanov}. The gauge-fixed forward action is bounded by the total
reference rate whatever the divergence, whereas the gauge-fixed reverse action
$\operatorname{KL}(r\|p)$ is unbounded; the two kernels measure two directions of
divergence, not two scales of the same one. We use the closed form of~(c) to
audit the same trained models at the forward kernel
(Section~\ref{sec:ex-audit}); the training penalty uses $\Psi$.
\end{remark}

\begin{remark}[What the penalty is]\label{rem:shrink}
Theorem~\ref{thm:vertex}(c) says that, state by state, the penalised minimiser is
a \emph{convex shrinkage} of the data posterior towards the amino-acid reference
of the language model, with weight $w$; the restriction to $V$ is what the sampler
uses (non-amino-acid tokens are excluded). For $\lambda=0.1$ and
$\nu_V\simeq0.93$, $w\simeq8.5\%$. The penalty therefore acts as an explicit,
measurable regulariser towards the protein-language-model prior ESM-2, of the form
of a mixture $(q+\lambda r)/(1+\lambda\nu_V)$; this is the sense in which it
lowers the action, and we do not interpret it further as a path-space optimality
statement. It is a population statement for an expressive class: a finite network
does not attain~\eqref{eq:shrink}. We make no claim that~\eqref{eq:gdsbm-loss} is
a path-space entropy: $\mathcal A_\Psi$ is a per-state functional, evaluated along
noised data states.
\end{remark}

\begin{remark}[Geometry, and what is not claimed]\label{rem:geometry}
The state space $\mathcal X=(V\cup\{\mathsf m\})^L$ with the edit edges is the
Cartesian product of $L$ stars $K_{1,|V|}$ centred at the mask; it has
$|\mathcal X|=(|V|+1)^L$ vertices and $|\mathcal E|=L|V|(|V|+1)^{L-1}$ edges. A
\emph{global} potential $v:\mathcal X\to\mathbb R$ yields an edge field
$v(y)-v(x)$ with zero circulation on every closed walk. Since each star is a
tree, the cycle space of the product is generated by the squares, so zero
circulation reduces to a commutation condition: filling position $i$ with $a$ then
position $j$ with $b$ must accumulate the same control as the reverse order. A
free edge field has $|\mathcal E|-|\mathcal X|+1$ additional rotational degrees of
freedom. The energy head, however, reads the potentials of the neighbours
\emph{from the current state}: consistency across states is learned, not
enforced, and we do not claim that $v_\theta$ is a global potential. What
Theorem~\ref{thm:vertex}(a) gives is the gauge pinning at each state, not the
curl-free condition.
\end{remark}

\begin{remark}[Scope of the conditioning claims]\label{rem:cond-scope}
Theorem~\ref{thm:cond} is a statement about the structure of the conditioned
model: a tilt is an $h$-transform, it stays in the vertex class, rewards compose,
the target of a statistic is closed-form, and a particle system reaches it. What
we measure in Section~\ref{sec:ex-cond} is the single-objective myopic sampler
against the closed-form marginal of part~(b); the particle system is stated but
not evaluated here, and neither is multi-objective guidance or the hydropathy
statistic, for which a reliable estimate of $\hat p$ would require far more than
$n_0=300$ unguided samples. In the particle system we would resample
systematically when the effective sample size falls below $K/4$ and report
weighted estimators without a final resampling step; we do \emph{not} use the
guided kernel~\eqref{eq:myopic-kernel} as proposal together with the
weights~\eqref{eq:incr-weight}, since the twist would then be applied twice, once
in the proposal and once in the weight. A guided proposal $q$ requires the
incremental weight to carry the ratio of the base kernel to $q$; we leave this
residual twist to future work and use the free-proposal version, whose target is
exactly~\eqref{eq:tempered} by Theorem~\ref{thm:cond}(c). Finally, the algebra
applies to both arms; the ordering observed in Section~\ref{sec:ex-cond} is
reported as measured and is not attributed to the vertex form rather than to the
penalty.
\end{remark}

\section{Girsanov's formula for jump processes}\label{app:girsanov}

We prove Fact~\ref{fact:girsanov} in full. Throughout, $\mathcal X$ is finite, the rates $R^0_t(x,y)$ are bounded measurable in $t$ and vanish off $\mathcal E$, and $u$ is bounded, so $R^u=R^0e^u$ is bounded and both $\mathbb P^0$ and $\mathbb P^u$ are laws of non-explosive Markov chains on $D([0,T],\mathcal X)$, the space of c\`adl\`ag trajectories from $[0,T]$ into $\mathcal X$.

\paragraph{Path likelihood.} Condition on $n$ jumps at times $t_1<\dots<t_n$ in $[0,T]$ with states $x_0\to x_1\to\dots\to x_n$. Writing $\Lambda_t(x)=\sum_{y\in \mathcal X}R_t(x,y)$ for the total rate, the likelihood of this path under rates $R$ with respect to the measure $\sum_{n\ge0}dt_1\cdots dt_n\,\delta$-counting on states is
\begin{equation}\label{eq:path-lik}
\mathcal L_R=\Big(\prod_{k=1}^{n}R_{t_k}(x_{k-1},x_k)\Big)\exp\Big(-\int_0^T\Lambda_t(X_t)\,dt\Big),
\end{equation}
the product of the jump intensities at the jump times and of the survival factors on the inter-jump intervals. The dominating measure does not depend on $R$, so the ratio of~\eqref{eq:path-lik} for $R^u$ and $R^0$ is well defined, and since $R^u$ and $R^0$ have the same support, $\mathbb P^u\ll\mathbb P^0$ with
\begin{equation}\label{eq:girs}
\begin{split}
\log\frac{d\mathbb P^u}{d\mathbb P^0}={}&\sum_{k=1}^nu_{t_k}(X_{t_k-},X_{t_k})\\
&-\int_0^T\sum_yR^0_t(X_t,y)\big(e^{u_t(X_t,y)}-1\big)\,dt,
\end{split}
\end{equation}
which is~\eqref{eq:girsanov-density}.

\paragraph{Compensation formula.} Let $N(dt,dy)$ be the jump measure of the chain, $N([0,t]\times\{y\})=\#\{k:t_k\le t,\ X_{t_k}=y\}$. Under a law with rates $R$, the compensator of $N$ is $R_t(X_t,y)\,dt$, i.e.\ for every bounded predictable $f$,
\begin{equation}\label{eq:compens}
\mathbb E\sum_{k}f(t_k,X_{t_k-},X_{t_k})=\mathbb E\int_0^T\sum_yR_t(X_t,y)\,f_t(X_t,y)\,dt .
\end{equation}
This is the definition of the intensity of a finite-state jump process and holds by Fubini on~\eqref{eq:path-lik}.

\paragraph{The two kernels.} Take $\mathbb E_{\mathbb P^u}$ of~\eqref{eq:girs}. The second term is already a time integral; for the first, apply~\eqref{eq:compens} under $R^u=R^0e^u$ with $f=u$:
\[
\operatorname{KL}(\mathbb P^u\|\mathbb P^0)=\mathbb E_{\mathbb P^u}\!\int_0^T\!\sum_yR^0\Big(u\,e^{u}-e^{u}+1\Big)dt=\mathbb E_{\mathbb P^u}\!\int_0^T\!\sum_yR^0\,a(u)\,dt,
\]
with $a(z)=ze^z-e^z+1$. Take instead $\mathbb E_{\mathbb P^0}$ of $-\log\frac{d\mathbb P^u}{d\mathbb P^0}$ and apply~\eqref{eq:compens} under $R^0$ with $f=u$:
\[
\operatorname{KL}(\mathbb P^0\|\mathbb P^u)=\mathbb E_{\mathbb P^0}\!\int_0^T\!\sum_yR^0\Big(-u+e^{u}-1\Big)dt=\mathbb E_{\mathbb P^0}\!\int_0^T\!\sum_yR^0\,\Psi(u)\,dt,
\]
with $\Psi(z)=e^z-z-1$. Both kernels are non-negative with a unique zero at $z=0$, and $\Psi(z)=e^{z}a(-z)$. Finally $\Psi''(z)=e^z>0$, so $\Psi$ is strictly convex, while $a''(z)=(z+1)e^z$ changes sign at $z=-1$, so $a$ is not convex. \qed

\section{The Gibbs--Doob Variational Principle}\label{app:gibbs}

This appendix proves Fact~\ref{fact:gibbs}. The result is classical; we give it in
full because the form of the optimal control, a difference of vertex values, is the
starting point of the whole paper, and because the cancellation that produces it is what
makes the $h$-transform a change of rates rather than a reweighting of paths.

\paragraph{Setting.} Keep the assumptions of Appendix~\ref{app:girsanov}: $\mathcal X$
finite, $R^0$ bounded on $[0,T]$, and now $\varphi:\mathcal X\to\mathbb R$ an arbitrary
function, necessarily bounded. Define
\[
h(x,t)=\mathbb E_{\mathbb P^0}\big[e^{\varphi(X_T)}\,\big|\,X_t=x\big],
\qquad v_t(x)=\log h(x,t),
\]
and let $\mathbb P^h$ be the law of the chain started at $x_0$ with rates
$R^0_t(x,y)h(y,t)/h(x,t)$, i.e.\ the controlled chain with
$u^\ast_t(x,y)=v_t(y)-v_t(x)$.

\paragraph{Step 0: $h$ is bounded away from zero.} Since $\mathcal X$ is finite,
$m=\min_x\varphi(x)$ and $M=\max_x\varphi(x)$ are finite, and $h(x,t)\in[e^{m},e^{M}]$ for
every $(x,t)$ as an average of values in that range. Hence $v=\log h$ is well defined and
bounded, $u^\ast$ is bounded, and $\mathbb P^h$ is the law of a non-explosive chain to
which Appendix~\ref{app:girsanov} applies.

\subsection*{Step 1: the backward equation}

\begin{lemma}\label{lem:backward}
$h$ is differentiable in $t$ and satisfies
\begin{equation}\label{eq:app-backward}
\partial_th(x,t)=-\sum_{y\in\mathcal X}R^0_t(x,y)\big(h(y,t)-h(x,t)\big),
\qquad h(\cdot,T)=e^{\varphi}.
\end{equation}
\end{lemma}

\begin{proof}
Condition on the behaviour of the chain on $[t,t+\delta]$. With probability
$1-\Lambda_t(x)\delta+o(\delta)$ there is no jump, and with probability
$R_t(x,y)\delta+o(\delta)$ there is exactly one jump, to $y$; two or more jumps have
probability $o(\delta)$ because the rates are bounded. The Markov property then gives
\[
h(x,t)=\big(1-\Lambda_t(x)\delta\big)h(x,t+\delta)
+\delta\sum_yR_t(x,y)h(y,t+\delta)+o(\delta).
\]
Subtracting $h(x,t+\delta)$, dividing by $\delta$ and letting $\delta\downarrow0$ yields
\eqref{eq:app-backward}; the terminal condition is the definition of $h$ at $t=T$.
\end{proof}

Dividing \eqref{eq:app-backward} by $h(x,t)>0$ and using
$u^\ast_t(x,y)=\log\frac{h(y,t)}{h(x,t)}$,
\begin{equation}\label{eq:app-dlogh}
\partial_tv_t(x)=-\sum_yR^0_t(x,y)\Big(\frac{h(y,t)}{h(x,t)}-1\Big)
=-\sum_yR^0_t(x,y)\big(e^{u^\ast_t(x,y)}-1\big).
\end{equation}

\subsection*{Step 2: telescoping along the path}

Fix a path with jumps at $t_1<\dots<t_n$. The function $t\mapsto v_t(X_t)$ is piecewise
$C^1$: between two jumps the state is constant and only $t$ varies, and at each jump time
it has a jump of size $v_{t_k}(X_{t_k})-v_{t_k}(X_{t_k-})=u^\ast_{t_k}(X_{t_k-},X_{t_k})$.
Summing its increments over $[0,T]$,
\[
v_T(X_T)-v_0(x_0)
=\sum_{k=1}^{n}u^\ast_{t_k}\big(X_{t_k-},X_{t_k}\big)
+\int_0^T\partial_tv_t(X_t)\,dt ,
\]
that is,
\begin{equation}\label{eq:app-telescope}
\sum_{k=1}^{n}u^\ast_{t_k}\big(X_{t_k-},X_{t_k}\big)
=\log h(X_T,T)-\log h(x_0,0)-\int_0^T\partial_tv_t(X_t)\,dt .
\end{equation}

\subsection*{Step 3: the cancellation}

Write the Girsanov density \eqref{eq:girsanov-density} for $u^\ast$:
\[
\log\frac{d\mathbb P^{h}}{d\mathbb P^0}
=\underbrace{\sum_{k}u^\ast_{t_k}}_{\text{Step 2}}
-\int_0^T\sum_yR^0_t(X_t,y)\big(e^{u^\ast_t(X_t,y)}-1\big)\,dt .
\]
Substituting \eqref{eq:app-telescope} for the jump term and
\eqref{eq:app-dlogh} for $\partial_tv_t(X_t)$,
\[
-\int_0^T\partial_tv_t(X_t)\,dt
=+\int_0^T\sum_yR^0_t(X_t,y)\big(e^{u^\ast_t(X_t,y)}-1\big)\,dt,
\]
which is exactly the compensator with the opposite sign. The two integrals cancel and
\begin{equation}\label{eq:app-doob-density}
\log\frac{d\mathbb P^{h}}{d\mathbb P^0}=\log h(X_T,T)-\log h(x_0,0)
=\varphi(X_T)-\log h(x_0,0),
\end{equation}
using the terminal condition $h(\cdot,T)=e^{\varphi}$. This is the first claim of
Fact~\ref{fact:gibbs}: the $h$-transform reweights a path by its terminal reward
alone, and the constant $h(x_0,0)$ is the normalisation.

\subsection*{Step 4: the variational identity}

Let $\mathbb P\ll\mathbb P^0$ be any path law started at $x_0$. On the event where the
densities are defined,
\[
\log\frac{d\mathbb P}{d\mathbb P^{h}}
=\log\frac{d\mathbb P}{d\mathbb P^{0}}-\log\frac{d\mathbb P^{h}}{d\mathbb P^{0}}
=\log\frac{d\mathbb P}{d\mathbb P^{0}}-\varphi(X_T)+\log h(x_0,0),
\]
by \eqref{eq:app-doob-density}. Taking $\mathbb E_{\mathbb P}$, and noting that
$\log h(x_0,0)$ is a constant,
\begin{equation}\label{eq:app-gibbs}
\operatorname{KL}\big(\mathbb P\,\|\,\mathbb P^{h}\big)
=\operatorname{KL}\big(\mathbb P\,\|\,\mathbb P^{0}\big)
-\mathbb E_{\mathbb P}\,\varphi(X_T)+\log h(x_0,0)\;\ge\;0 ,
\end{equation}
which is \eqref{eq:gibbs}. Rearranging,
\[
\operatorname{KL}\big(\mathbb P\|\mathbb P^{0}\big)-\mathbb E_{\mathbb P}\varphi(X_T)
=-\log h(x_0,0)+\operatorname{KL}\big(\mathbb P\|\mathbb P^{h}\big),
\]
so the left-hand side is bounded below by $-\log h(x_0,0)$, with equality if and only if
$\operatorname{KL}(\mathbb P\|\mathbb P^{h})=0$, i.e.\ $\mathbb P=\mathbb P^{h}$. The
minimiser is therefore unique and is the $h$-transform. \qed

\subsection*{Remark: why the optimal control is a gradient}

Nothing in Steps 1 to 4 used the structure of $\mathcal X$ beyond finiteness. What the
computation shows is that the optimal control is $u^\ast_t(x,y)=v_t(y)-v_t(x)$ with $v$ a
function of the \emph{vertex} alone, and that this is forced: it is the only form for
which the telescoping of Step 2 can cancel the compensator of Step 3. An edge field that
is not a difference of vertex values leaves a residue in \eqref{eq:app-doob-density}, and
the resulting law is not an $h$-transform of $\mathbb P^0$. This is the observation that
Section~\ref{sec:gdsbm} turns into a parameterisation.

\section{Proof of Theorem~\ref{thm:gauge}}\label{app:gauge}

Everything in this appendix is one identity read three times. The per-state action splits
into a divergence that the conditional determines and a term that only the gauge sees;
proving this once gives, in turn, the value of the gauge that minimises the action and the
corresponding closed form for the forward Girsanov kernel.

\paragraph{Notation.} Fix $(x,t)$ and write $\nu=\sum_yR^0(x,y)$, $r=R^0(x,\cdot)/\nu$,
$Z=\sum_yR^0e^{u}$, $\zeta=Z/\nu=\sum_yr(y)e^{u(y)}$, $p=R^0e^{u}/Z=re^{u}/\zeta$, and
$h(\zeta)=\zeta-1-\log\zeta$ for $\zeta>0$. Recall $\Psi(z)=e^{z}-z-1$ and $a(z)=ze^{z}-e^{z}+1$.

\subsection*{Part (a): the splitting}

\begin{proof}
Everything follows from writing the control in terms of the conditional it produces.
From $p=R^0e^{u}/Z$ we get
\[
u=\log\frac{p}{R^0}+\log Z ,
\]
and the three sums that make up $\Psi(u)=e^{u}-u-1$ can then be computed one by one.
The exponential sum is the normalisation itself, $\sum_yR^0e^{u}=Z=\nu\zeta$. The constant
sum is $\sum_yR^0=\nu$. For the linear sum, substitute the expression for $u$ and factor
out $\nu$:
\[
\sum_{y}R^0\,u
=\sum_{y}R^0\log\frac{p}{R^0}+\nu\log Z
=\nu\sum_{y}r\,\log\frac{p}{\nu r}+\nu\log Z
=-\nu\operatorname{KL}(r\|p)-\nu\log\nu+\nu\log Z ,
\]
where the middle equality uses $R^0=\nu r$ and the last one separates the constant
$\log\nu$ from the divergence. Assembling the three pieces,
\[
\sum_{y}R^0\,\Psi(u)
=\nu\zeta-\nu+\nu\operatorname{KL}(r\|p)+\nu\log\nu-\nu\log Z
=\nu\Big[\operatorname{KL}(r\|p)+\zeta-1-\log\zeta\Big],
\]
since $\log Z-\log\nu=\log\zeta$. This is \eqref{eq:gauge-decomp}. Non-negativity of the
first term is that of a relative entropy. For the second, the inequality
$\log\zeta\le\zeta-1$ gives $h\ge0$ on $(0,\infty)$, with equality only at $\zeta=1$;
moreover $h(\zeta)\to\infty$ both as $\zeta\downarrow0$ and as $\zeta\to\infty$. The first
term is a function of $p$ alone and the second of $\zeta$ alone.
\end{proof}

\subsection*{Part (b): the orbit}

\begin{proof}
The cross-entropy~\eqref{eq:ce} depends on $u$ only through $p$. Adding a constant $c$ to
$u$ multiplies both $R^0e^{u}$ and $Z$ by $e^{c}$, so $p$ is unchanged while $\zeta$
becomes $\zeta e^{c}$; as $c$ runs over $\mathbb R$, $\zeta e^{c}$ runs over $(0,\infty)$
and $h(\zeta e^{c})$ over $[0,\infty)$. By part~(a) the action along the orbit is
$\nu\operatorname{KL}(r\|p)+\nu h(\zeta e^{c})$, whose minimum over the gauge is therefore
reached exactly when $\zeta e^{c}=1$, that is at $c^\ast=-\log\zeta$, with value
$\nu\operatorname{KL}(r\|p)$, and whose excess over that minimum is $\nu h(\zeta e^{c})$,
unbounded in $c$.
\end{proof}

\subsection*{Part (c): the forward kernel}

\begin{proof}
It suffices to treat $\nu=1$ and multiply by $\nu$ at the end, since $R^0=\nu r$.
From $e^{u}=e^{c}p/r$ we get $\sum_yr\,e^{u}=e^{c}$ and
$\sum_yr\,u\,e^{u}=e^{c}\sum_yp\big(\log(p/r)+c\big)=e^{c}\big(\operatorname{KL}(p\|r)+c\big)$,
so $\sum_yr\,a(u)=e^{c}\operatorname{KL}(p\|r)+ce^{c}-e^{c}+1$. Its derivative in $c$ is
$e^{c}\big(\operatorname{KL}(p\|r)+c\big)$, which vanishes and changes sign at
$c^\ast=-\operatorname{KL}(p\|r)$; substituting gives $1-e^{-\operatorname{KL}(p\|r)}$.
The map $d\mapsto1-e^{-d}$ is strictly increasing from $[0,\infty)$ onto $[0,1)$, hence a
bijection, so a measured gauge-fixed forward action recovers
$\operatorname{KL}(p\|r)=-\log(1-\cdot)$.
\end{proof}

\section{Proof of Theorem~\ref{thm:vertex}}\label{app:vertex}

\subsection*{Part (a): pinning the gauge costs nothing}

\begin{proof}
By definition $\log\big(p_\theta[u](a)/r(a)\big)=u_a-\log Z_u$ for every $a$. If
$u-w=c\mathbf 1$, then $Z_u=e^{c}Z_w$ and the two conditionals coincide. Conversely, if
$p_\theta[u]=p_\theta[w]$, then $u_a-w_a=\log Z_u-\log Z_w$ for every $a$, a constant
independent of $a$, so $u-w\in\mathbb R\mathbf 1$. For $c=-u_s$ this gives
$u-u_s\mathbf 1\in\mathcal U_{\rm vtx}$ and
\[
p_\theta[u-u_s\mathbf 1](a)=\frac{r(a)\,e^{u_a-u_s}}{\sum_b r(b)\,e^{u_b-u_s}}=\frac{e^{-u_s}\,r(a)\,e^{u_a}}{e^{-u_s}\,Z_u}=p_\theta[u](a),
\]
so $p_\theta$ is the same. Uniqueness holds because two elements of
$\mathcal U_{\rm vtx}$ that differ by $c\mathbf 1$ have $s$-components $0$ and $0+c$, hence
$c=0$. If $u_a=v_a-v_s$, then $u_s=0$; conversely, if $u_s=0$, take $v=u$. Here $v$ is a
potential on the star of neighbours of the current state, not a global potential on
$\mathcal X$ (Remark~\ref{rem:geometry}). Finally, inverting
$p_\theta[u](a)=r(a)e^{u_a}/Z_u$ in the vertex gauge, where $u_s=0$ gives
$Z_u=r(s)/p_s$, yields $u_a=\log\big[(p_a/r(a))/(p_s/r(s))\big]$.

The map $u\mapsto p_\theta[u]$ is therefore constant on each orbit $u+\mathbb R\mathbf 1$,
and every orbit meets $\mathcal U_{\rm vtx}$ exactly once, so the two classes produce the
same set of conditionals. Since $\mathrm{CE}$ depends on $u$ only through $p_\theta[u]$,
its infimum is the same over both classes, and GDSBM at $\lambda=0$ is functionally
identical to MadSBM. What the conditionals do not fix is the rate level: $R^0e^{u}$ and
$R^0e^{u-u_s\mathbf 1}$ differ by the factor $e^{-u_s}$, which is the part of the action
that Theorem~\ref{thm:gauge} isolates.
\end{proof}

\subsection*{Part (b): the anchor}

\begin{proof}
Theorem~\ref{thm:gauge}(a) with $\nu=\sum_yr(y)=1$ gives
$\mathcal A_\Psi(u;x)=\operatorname{KL}(r\|p)+h(\zeta)$ with
$\zeta=\sum_yr(y)e^{u(y)}$ and $p=re^u/\zeta$. In the vertex gauge $u(s)=0$, so
$p(s)=r(s)e^{u(s)}/\zeta=r(s)/\zeta$, i.e.\ $\zeta=r(s)/p(s)$, which is
\eqref{eq:vgauge}. Since $h$ vanishes only at $1$, the second term vanishes if and only if
$p(s)=r(s)$.
\end{proof}

\subsection*{Part (c): the penalised per-state minimiser}

Fix a masked position with current token $s=\mathsf m$, let $r>0$ on $\mathcal T$, let $q$ be the clean-token posterior, supported on $V$ so that $q(s)=0$, and let $\nu=1$ (unit intensity). Minimise over the open simplex
\[
J(p)=-\sum_{y}q(y)\log p(y)+\lambda\Big[\operatorname{KL}(r\|p)+h\big(r(s)/p(s)\big)\Big].
\]

\paragraph{Step 1: an explicit form.} Using $\operatorname{KL}(r\|p)=\sum_yr\log r-\sum_yr\log p$ and $h(r(s)/p(s))=r(s)/p(s)-1-\log r(s)+\log p(s)$,
\begin{equation}\label{eq:app-J}
J(p)=-\sum_{y\neq s}\big(q(y)+\lambda r(y)\big)\log p(y)\;+\;\lambda\big(1-r(s)\big)\log p(s)\;+\;\lambda\,\frac{r(s)}{p(s)}\;+\;\mathrm{const},
\end{equation}
where the constant collects $\sum_yr\log r$, $-1$ and $-\log r(s)$, and we used $q(s)=0$.

\paragraph{Step 2: existence of an interior minimiser.} $J$ is smooth on the open simplex. If $p(y)\to0$ for some $y\neq s$, then $-\log p(y)\to+\infty$ with coefficient $q(y)+\lambda r(y)\ge\lambda r(y)>0$, so $J\to+\infty$. If $p(s)\to0$, then $\lambda r(s)/p(s)\to+\infty$ like $1/p(s)$ while $\lambda(1-r(s))\log p(s)\to-\infty$ only like $\log p(s)$, so again $J\to+\infty$. Hence the sublevel sets are compact in the open simplex and a minimiser exists in the interior.

\paragraph{Step 3: stationarity.} At an interior minimiser there is $\kappa>0$ with $\partial J/\partial p(y)=-\kappa$ for all $y$. For $y\neq s$ this reads $-\big(q(y)+\lambda r(y)\big)/p(y)=-\kappa$, i.e.
\begin{equation}\label{eq:app-stat}
p(y)=\frac{q(y)+\lambda r(y)}{\kappa},\qquad y\neq s,
\end{equation}
while for $y=s$ it gives the quadratic $\kappa\,p(s)^2+\lambda\big(1-r(s)\big)p(s)-\lambda r(s)=0$, which has a unique positive root; $\kappa$ is then fixed by $\sum_yp(y)=1$.

\paragraph{Step 4: restriction to the amino acids.} The sampler emits only $a\in V$, and $V\subset\mathcal T\setminus\{s\}$, so~\eqref{eq:app-stat} applies to all of $V$. Summing over $V$ and using $\sum_{a\in V}q(a)=1$ gives $\sum_{a\in V}p(a)=(1+\lambda\nu_V)/\kappa$ with $\nu_V=\sum_{a\in V}r(a)$. Dividing, $\kappa$ cancels:
\begin{equation}\label{eq:app-shrink}
p(a\mid V)=\frac{q(a)+\lambda\nu_V\tilde r(a)}{1+\lambda\nu_V}=(1-w)\,q(a)+w\,\tilde r(a),\qquad w=\frac{\lambda\nu_V}{1+\lambda\nu_V},\quad\tilde r=\frac{r|_V}{\nu_V}.
\end{equation}
The total-variation statement follows from $p(\cdot\mid V)-q=w(\tilde r-q)$. \qed

\paragraph{Numerical value.} The audit of Section~\ref{sec:ex-audit} measures $\nu_V\in[0.933,0.936]$ along the sampling trajectories, so at $\lambda=0.1$ the shrinkage weight is $w=0.085$. Equation~\eqref{eq:app-shrink} is a statement about the population minimiser over an expressive class; a finite network does not attain it, and we do not use it as a prediction of the trained conditional.

\section{Proof of Theorem~\ref{thm:cond}}\label{app:cond}

\subsection*{Part (a): the tilt is an $h$-transform, and composition}

\begin{proof}
Apply Fact~\ref{fact:gibbs} with $\mathbb P^\ast$ in place of $\mathbb P^0$ and
$\varphi=\beta\Phi$: $d\mathbb P^\beta/d\mathbb P^\ast=e^{\beta\Phi(X_T)}/h^{\beta\Phi}(x_0,0)$,
so the terminal marginal of $\mathbb P^\beta$ is $\mu_1e^{\beta\Phi}/Z_\beta=\mu_1^\beta$
(and $h^{\beta\Phi}(x_0,0)=Z_\beta$ for a Dirac source). The rates of $\mathbb P^\beta$
relative to $\mathbb P^0$ are $R^0e^{u^\ast}e^{v^{\beta\Phi}(y)-v^{\beta\Phi}(x)}$, which
is~\eqref{eq:add-control}; a sum of vertex potentials is a vertex potential.

For the composition, the densities multiply:
$e^{\Phi_2}\,d\mathbb P^{\Phi_1}\propto e^{\Phi_1+\Phi_2}\,d\mathbb P^\ast$, so tilting by
$\Phi_1$ then by $\Phi_2$ gives $\mathbb P^{\Phi_1+\Phi_2}$. Conditionally on $X_t=x$,
$d\mathbb P^{\Phi_1}/d\mathbb P^{\ast}=e^{\Phi_1(X_T)}/\mathbb E_{\mathbb P^\ast}[e^{\Phi_1}\mid X_t=x]$,
hence
$\mathbb E_{\mathbb P^\ast}[e^{\Phi_1+\Phi_2}\mid x]=\mathbb E_{\mathbb P^\ast}[e^{\Phi_1}\mid x]\cdot\mathbb E_{\mathbb P^{\Phi_1}}[e^{\Phi_2}\mid x]$;
taking logarithms gives~\eqref{eq:compose-gap}. At $t=T$ both conditional expectations
equal $e^{\Phi_2(x)}$, so $G(\cdot,T)=0$.
\end{proof}

\subsection*{Part (b): the closed-form marginal}

\begin{proof}
$\mu_1^\beta(c=v)=Z_\beta^{-1}\sum_{x:\,c(x)=v}\mu_1(x)\,e^{\beta\psi(c(x))}=Z_\beta^{-1}e^{\beta\psi(v)}\,\mu_1(c=v)$,
and $Z_\beta=\sum_{v'}e^{\beta\psi(v')}\mu_1(c=v')$ by the same grouping. Only the fact
that $\Phi$ is constant on the level sets of $c$ is used.
\end{proof}

\subsection*{Part (c): bias and telescoping weights}

\begin{proof}
The kernel~\eqref{eq:myopic-kernel} reweights the one-step conditional by
$e^{\beta\Delta\Phi}$, which is the exact tilt of that conditional and not of the terminal
law, since the exact terminal tilt would require $v^{\beta\Phi}$ and not $\beta\Phi$; the
discrepancy is the gap of part~(a) at $t<T$ and is the bias documented in
\citep{lee2025debiasing,wan2026dgm}.

For the particle system, the product~\eqref{eq:incr-weight} telescopes to
$e^{\beta[\Phi(x^{(N)})-\Phi(x^{(0)})]}$, and the constant $e^{-\beta\Phi(x^{(0)})}$
cancels under self-normalisation; this uses that $x^{(N)}$ is the fully decoded sequence,
so that $\Phi(x^{(N)})$ is the reward of the returned sample. The weighted sample is then
self-normalised importance sampling from the proposal $\mu_1$ with weight
$e^{\beta\Phi}$, whose normalised target is~\eqref{eq:tempered}; consistency is the
standard law of large numbers for self-normalised estimators ($\Phi$ is bounded), and
resampling according to the current weights leaves the Feynman--Kac flow invariant
\citep{wu2023tds,skreta2025feynman}.
\end{proof}

\section{Experimental details, deviations and threats to validity}\label{app:details}

\paragraph{Architecture.} Frozen ESM-2 650M \citep{lin2023esm} provides logits $f_\phi(x_t)$ and last-layer hidden states. The trainable head is two AdaLN transformer blocks (hidden $1280$, $16$ heads, MLP ratio $4$, $\approx45$M parameters) conditioned on a Gaussian Fourier time embedding ($64$ dimensions, scale $30$) followed by a SiLU MLP, then a final norm and a vocabulary projection whose weight and bias are initialised to zero. At initialisation the control vanishes and the model is exactly the tempered reference. The MadSBM arm reads the projection as a free field $u_\theta(x,\cdot,t)$; the GDSBM arm reads it as a vertex potential $v_\theta(x,\cdot,t)$ and forms $u_\theta=\bar\nabla v_\theta$ as in Definition~\ref{def:vertex}.

\paragraph{Optimisation.} AdamW, $(\beta_1,\beta_2)=(0.9,0.999)$, weight decay $0.01$, learning rate $10^{-4}$ with cosine decay to $10^{-5}$, warmup over $4\%$ of steps, gradient clipping at $10$, batch size $64$, $127\,000$ steps, $500$ discrete masking levels with $t$ uniform and clamped to $[10^{-6},1-10^{-6}]$, validation every $500$ steps with checkpoint selection on the logged objective. $\lambda=0.1$ for GDSBM and $\lambda=0$ for MadSBM.

\paragraph{Sampling.} MadSBM's sampler, unchanged for both arms: rate scale $\kappa=0.01$, jump scale $\eta=0.05$, temperature $\tau=0.5$, top-$p$ $\pi_{\rm top}=0.9$, rate clamp $-40$, $N$ steps from $t_{\max}=1-10^{-6}$ down to $10^{-6}$, with positions still masked at the end resolved by the final argmax. Non-amino-acid tokens are excluded from the candidate set; a row emptied by top-$p$ filtering is recomputed without the filter, which occurred for about one position in the whole campaign.

\begin{algorithm}[t]
\caption{Sampling (identical for both arms)}
\label{alg:sample}
\begin{algorithmic}[1]
\Require trained head $\theta$, length $L$, budget $N$; $t$ decreases from $t_{\max}$ to $\epsilon$, $\mathrm dt=1/N$
\State $x\gets\mathsf m^{L}$
\For{$k=1,\dots,N$}
  \State read $u_i$ from the head as in line~\ref{line:arm} of Algorithm~\ref{alg:train}; $\ell_i\gets\beta_tf_\phi(x)_i+u_i$
  \State $\pi_i\gets1-\exp\big(-\min\{\eta\,\textstyle\sum_{y}e^{\kappa u_i(y)}\,\mathrm dt,\,40\}\big)$ \label{line:intensity}
  \State $a_i\sim\operatorname{softmax}\big(\mathrm{top}\text{-}p_{\pi_{\rm top}}(\ell_i/\tau)\big|_V\big)$ \label{line:kernel}
  \State $x_i\gets a_i$ for each masked $i$ that jumps, with probability $\pi_i$
\EndFor
\State fill the positions still masked by $\arg\max_{y\in V}\ell_i(y)$; \Return $x$
\end{algorithmic}
\end{algorithm}

\paragraph{Metrics.} Pseudo-perplexity under ESM-2 650M in fp32 by sequential single-position masking. pLDDT from \texttt{EsmForProteinFolding} in fp32, one sequence at a time, averaged per sequence and rescaled to $[0,100]$. Uniqueness is the fraction of distinct sequences in a generation set. Action curves are computed along the sampling trajectories with the closed forms of Appendix~\ref{app:gauge}, averaged over internal positions, on $46$ sequences per (arm, seed, $N$).

\paragraph{Declared deviations from \citep{goel2026madsbm}.} (i)~Corpus rebuilt from UniProt (reviewed, $5\le L\le50$, strict 20-letter alphabet), their dataset not being public; comparisons are internal to our two arms. (ii)~Standard padding with batch $64$ instead of dynamic batching at batch $1$ with accumulation $16$. (iii)~Our own training loop. (iv)~Conditional guidance removed from the sampler for the unconditional comparison. (v)~pLDDT through \texttt{transformers} rather than OpenFold, same model and precision. (vi)~Two sampler safeguards are ours and are disclosed: non-amino-acid tokens are excluded from the candidate set, and a row emptied by top-$p$ filtering is recomputed without the filter. (vii)~The conditioning runs of Section~\ref{sec:ex-cond} were executed on CPU. Their engineering constants were otherwise kept unchanged for comparability.

\paragraph{Checkpoint selection and excluded run.} The \texttt{madsbm} run at seed $44$ was resumed twice and its best-step tracker was reset, so the retained checkpoint is at step $124\,000$ rather than at the best step $76\,000$; this is the one selection asymmetry in the campaign. A $300$-step diagnostic run (\texttt{madsbm} at seed $97$) is present in our logs and is excluded from every table, as a run of a different length from every reported configuration.

\paragraph{Statistical scope.} Three training seeds, hence two degrees of freedom in every paired test. We report constancy of sign, uncorrected $p$-values and the Bonferroni threshold, and we do not claim significance for any single metric of Section~\ref{sec:ex-uncond}. The claims we do make are the identity of Theorem~\ref{thm:vertex}(a) verified to float precision, the gauge decomposition of Theorem~\ref{thm:gauge}(a) verified numerically through the agreement of the gauge-fixed action with $\operatorname{KL}(r\|p)$, and the ordering of the gauge-fixed actions, which holds in all nine (seed, $N$) pairs at both kernels. The conditioning measurement is on three seeds, one checkpoint per seed and arm, executed on CPU.

\paragraph{What the action reduction is not.} The penalty uses the reverse kernel and is evaluated per state along noised data; Theorem~\ref{thm:vertex}(c) shows that it acts as a shrinkage towards the reference with weight $w\simeq8.5\%$ at $\lambda=0.1$. A lower gauge-fixed action therefore means a conditional closer to the tempered language model in the reverse-KL sense, not that the trained path measure minimises the path-space entropy of the Schr\"odinger problem.

\section{Numerical tables}\label{app:tables}

This appendix gives the numerical values behind Figure~\ref{fig:action-gauge} and behind Section~\ref{sec:ex-uncond}; the conditioning values are in Table~\ref{tab:cond} of the main text.

\begin{table}[h]
\centering\footnotesize
\caption{\textbf{Action audit along sampling trajectories}, mean $\pm$ s.d.\ over three seeds. $\mathcal A_{\rm raw}$ is the action of the control as parameterised; $\mathcal A_{\rm vtx}$ is the same control put in the vertex gauge $u(x_i)=0$; $\mathcal A_{\rm gf}$ is the minimum over the gauge class, which equals $\operatorname{KL}(r\|p)$ by Theorem~\ref{thm:gauge}. Left block: full $33$-token vocabulary. Right block: restricted to the $20$ amino acids, the only tokens the sampler can emit. For GDSBM, $\mathcal A_{\rm raw}=\mathcal A_{\rm vtx}$ by construction. The raw columns are seed means; their inter-seed dispersion is reported in the text.}
\label{tab:gauge}
\resizebox{\textwidth}{!}{%
\begin{tabular}{@{}llccc@{\hskip 1.2em}ccc@{}}
\toprule
& & \multicolumn{3}{c}{33 tokens} & \multicolumn{3}{c}{20 amino acids}\\
\cmidrule(r){3-5}\cmidrule(l){6-8}
Model & $N$ & $\mathcal A_{\rm raw}$ & $\mathcal A_{\rm vtx}$ & $\mathcal A_{\rm gf}$ & $\mathcal A_{\rm raw}$ & $\mathcal A_{\rm vtx}$ & $\mathcal A_{\rm gf}$\\
\midrule
MadSBM & 32 & $2.9\!\cdot\!10^{8}$ & $5.7\!\cdot\!10^{8}$ & $6.35\pm0.92$ & $2.9\!\cdot\!10^{8}$ & $5.7\!\cdot\!10^{8}$ & $5.25\pm0.83$\\
MadSBM & 64 & $3.5\!\cdot\!10^{8}$ & $6.1\!\cdot\!10^{8}$ & $6.51\pm1.15$ & $3.5\!\cdot\!10^{8}$ & $6.1\!\cdot\!10^{8}$ & $5.38\pm1.08$\\
MadSBM & 128 & $8.9\!\cdot\!10^{8}$ & $9.1\!\cdot\!10^{8}$ & $6.72\pm1.31$ & $8.9\!\cdot\!10^{8}$ & $9.1\!\cdot\!10^{8}$ & $5.58\pm1.23$\\
\midrule
GDSBM & 32 & $1.808\pm0.039$ & $1.808\pm0.039$ & $1.512\pm0.040$ & $1.655\pm0.034$ & $1.655\pm0.034$ & $1.369\pm0.034$\\
GDSBM & 64 & $1.826\pm0.117$ & $1.826\pm0.117$ & $1.514\pm0.096$ & $1.670\pm0.117$ & $1.670\pm0.117$ & $1.370\pm0.097$\\
GDSBM & 128 & $1.852\pm0.117$ & $1.852\pm0.117$ & $1.542\pm0.097$ & $1.696\pm0.114$ & $1.696\pm0.114$ & $1.399\pm0.095$\\
\bottomrule
\end{tabular}}
\end{table}

\begin{table}[h]
\centering\small
\caption{\textbf{Forward-kernel audit}, gauge-fixed, mean $\pm$ s.d.\ over three seeds, $20$ amino acids. $\mathcal A^{a}_{\rm gf}=1-e^{-\operatorname{KL}(p\|r)}$ by Theorem~\ref{thm:gauge}(c) at unit intensity; the last column inverts it on the seed means. The ordering holds in all nine (seed, $N$) pairs.}
\label{tab:fwd}
\resizebox{\textwidth}{!}{%
\begin{tabular}{@{}lcccc@{}}
\toprule
$N$ & MadSBM $\mathcal A^{a}_{\rm gf}$ & GDSBM $\mathcal A^{a}_{\rm gf}$ & ratio & ratio of $\operatorname{KL}(p\|r)$\\
\midrule
32  & $0.5494\pm0.0091$ & $0.5094\pm0.0072$ & $1.08$ & $1.12$\\
64  & $0.5503\pm0.0164$ & $0.5125\pm0.0207$ & $1.07$ & $1.11$\\
128 & $0.5613\pm0.0217$ & $0.5170\pm0.0183$ & $1.09$ & $1.13$\\
\bottomrule
\end{tabular}}
\end{table}

\begin{table}[h]
\centering\small
\caption{\textbf{Unconditional generation}, $920$ sequences per (arm, seed, $N$); mean over seeds of the per-seed means $\pm$ inter-seed s.d.}
\label{tab:uncond}
\resizebox{\textwidth}{!}{%
\begin{tabular}{@{}llcccc@{}}
\toprule
$N$ & Model & seeds & PPL $\downarrow$ & pLDDT $\uparrow$ & unique $\uparrow$\\
\midrule
32 & MadSBM (reimpl.) & 3 & $9.550\pm0.256$ & $73.44\pm0.69$ & $0.919\pm0.016$\\
32 & GDSBM (ours) & 3 & $9.796\pm0.025$ & $72.95\pm0.17$ & $0.928\pm0.021$\\
64 & MadSBM (reimpl.) & 3 & $9.520\pm0.091$ & $73.37\pm0.58$ & $0.909\pm0.029$\\
64 & GDSBM (ours) & 3 & $9.475\pm0.226$ & $73.13\pm0.60$ & $0.919\pm0.027$\\
128 & MadSBM (reimpl.) & 3 & $9.488\pm0.274$ & $73.64\pm0.40$ & $0.914\pm0.032$\\
128 & GDSBM (ours) & 3 & $9.529\pm0.035$ & $73.52\pm0.20$ & $0.912\pm0.027$\\
\bottomrule
\end{tabular}}
\end{table}

\begin{table}[h]
\centering\small
\caption{\textbf{Re-evaluated cross-entropy}, common random numbers, three seeds. Paired at equal seed, the difference GDSBM${}-{}$MadSBM on the test split is $+0.040$, $-0.125$ and $-0.053$ for seeds $42,43,44$: the sign is not constant, so the penalty neither costs nor buys cross-entropy systematically at $\lambda=0.1$.}
\label{tab:stability}
\begin{tabular}{@{}llccc@{}}
\toprule
Split & Model & seeds & CE & CV(CE)\\
\midrule
val & MadSBM (reimpl.) & 3 & $1.1825\pm0.0910$ & $0.077$\\
val & GDSBM (ours) & 3 & $1.1367\pm0.0449$ & $0.039$\\
test & MadSBM (reimpl.) & 3 & $1.1915\pm0.0741$ & $0.062$\\
test & GDSBM (ours) & 3 & $1.1454\pm0.0308$ & $0.027$\\
\bottomrule
\end{tabular}
\end{table}

\section{Deferred theory}\label{app:future}

Three questions are left open; we state them here because they bound the present contribution rather than extend it.

\paragraph{From a per-state regulariser to a path-space action.} Our penalty is the reverse Girsanov kernel evaluated state by state along noised data, not $\operatorname{KL}(\mathbb P^0\|\mathbb P^u)$ along the controlled process. Two routes could close the gap: evaluate the kernel along sampled trajectories, at the cost of differentiating through the sampler; or penalise the forward kernel $a$ with the gauge fixed in closed form (Theorem~\ref{thm:gauge}(c)), which is the correct path-space integrand but whose kernel is not convex, so the optimisation needs study.

\paragraph{From a locally read potential to a global one.} The energy head reads the potentials of the neighbours from the current state, so nothing forces the learned field to be curl-free on the squares of the edit graph (Remark~\ref{rem:geometry}). A penalty on the circulation around squares is computable in closed form on a minibatch and would make the control an exact discrete gradient; whether this changes the generated distribution is an empirical question we have not settled. Related to it, a bridge whose reference is the trained conditional rather than the frozen language model would make the whole construction self-consistent, at the cost of an outer iteration.

\paragraph{Geometry of the control potential.} We read the potential locally, as a vector on the neighbours of the current state, and make no use of the geometry of the sequence graph. A natural next step is to equip the edit graph with the Onsager-type structure that governs entropy dissipation for finite-state Schr\"odinger bridges \citep{benabdallah2026nonlocal}, and to ask whether a potential read in that geometry, rather than by a state-wise head, yields a control that is globally consistent and whose action carries a dissipation interpretation. This would address the two questions above at once.

\end{document}